\documentclass[a4paper, 11pt]{amsart} 
\usepackage[
  margin=1.9cm,
  includefoot,
  footskip=30pt,
]{geometry}
\usepackage{amsmath,amssymb,amscd}
\usepackage{amsfonts}
\usepackage[english]{babel}
\usepackage[leqno]{amsmath}
\usepackage{amssymb,amsthm}
\usepackage{mathrsfs} 
\usepackage{amscd}
\usepackage{enumerate}
\usepackage{epsfig}
\usepackage{color}
\usepackage{relsize}
\usepackage{layout}
\usepackage[usenames,dvipsnames]{xcolor}
\usepackage[backref=page]{hyperref}
\usepackage{tikz}
\usepackage{tikz-cd}
\hypersetup{
 colorlinks,
 citecolor=Green,
 linkcolor=Red,
 urlcolor=Blue}
\usepackage[matrix,arrow,tips,curve]{xy}
\input{xy}
\xyoption{all}
\definecolor{darkgreen}{rgb}{0.0, 0.7, 0.0}
\definecolor{purple}{rgb}{0.5, 0.0, 0.5}
\definecolor{red}{rgb}{0.8, 0.2, 0.0}

\newtheorem{thm}{Theorem}[section]
\newtheorem{bthm}{Theorem}

\newtheorem{lemma}[thm]{Lemma}

\newtheorem{prop}[thm]{Proposition}

\newtheorem{claim}[thm]{Claim}

\numberwithin{equation}{section}
\theoremstyle{definition}

\newtheorem{defi-not}[thm]{Definition-Notation}

\newtheorem{notation}[thm]{Notation}

\theoremstyle{remark}
\newtheorem{remark}[thm]{Remark}

\newtheorem{example}[thm]{Example}
\newcommand{\Z}{\mathbb{Z}}
\newcommand{\C}{\mathbb{C}}
\newcommand{\Q}{\mathbb{Q}}

\newcommand{\Pic}{\operatorname{Pic}}

\DeclareMathOperator{\Ext}{{Ext}}

\def \Im{{\rm Im}}
\DeclareMathOperator{\Supp}{Supp}

\def \P{\mathbb{P}}

\def \F{\mathcal F}
\def \N {\mathcal N}

\def \K{\mathcal K}

\def\I{{\mathcal J}}
\def \L{\mathcal L}
\def \E{\mathcal E}
\def \G{\mathcal G}

\def\O{\mathcal O}

\def \V{\mathcal V}

\def\M0{\mathcal M^0}

\def\mapright#1{\smash{\mathop{\longrightarrow}\limits^{#1}}}

\newcommand{\Id}{{\operatorname{Id}}}

\DeclareMathOperator{\codim}{{codim}}

\DeclareMathOperator{\Sing}{{Sing}}

\DeclareMathOperator{\Ker}{{Ker}}

\DeclareMathOperator{\Proj}{{Proj}}
\DeclareMathOperator{\Sym}{{Sym}}

\newcommand{\rk}{\operatorname{rank}}
\newcommand{\ch}{\operatorname{ch}}

\newcommand{\Pseff}{\operatorname{Pseff}}
\newcommand{\Nef}{\operatorname{Nef}}

\begin{document} 

\title[On the classification of non-big Ulrich bundles on varieties with Picard number one]{On the classification of non-big Ulrich bundles on varieties with Picard number one}

\author[A.F. Lopez, Y. Takahashi]{Angelo Felice Lopez and Yuta Takahashi}

\address{\hskip -.43cm Angelo Felice Lopez, Dipartimento di Matematica e Fisica, Universit\`a di Roma
Tre, Largo San Leonardo Murialdo 1, 00146, Roma, Italy. e-mail {\tt angelo.lopez@uniroma3.it}}

\address{\hskip -.43cm Yuta Takahashi, Department of Mathematics, Faculty of Science and Engineering, Chuo University. 1-13-27 Kasuga, Bunkyo-ku, Tokyo 112-8551, Japan. email: {\tt yuta0630takahashi0302@gmail.com}}

\thanks{The first author is partially supported by the GNSAGA group of INdAM. The second author is partially supported by JST SPRING, Japan Grant Number JPMJSP2170.}


\begin{abstract}
We characterize the triples $(X,\O_X(1),\E)$, where $X$ is a smooth projective variety with Picard number one and $\E$ is a non-big Ulrich bundle on $X$.
\end{abstract}

\maketitle

\section{Introduction}
Ulrich bundles are important objects in the study of projective varieties because of their strong cohomological properties and their connections with linear resolutions. In particular, their existence, construction and classification have been widely studied. From the viewpoint of positivity, recall that every Ulrich bundle is globally generated. Hence the tautological line bundle $\O_{\P(\E)}(1)$ is nef, and it is natural to ask whether it is also big. The positivity of Ulrich bundles was studied in \cite{lo, lm, lms1, lms2}. The results in these papers suggest that the failure of bigness imposes strong constraints on the geometry of the underlying variety.

Let $X$ be a smooth irreducible projective variety of dimension $n$. As is well known, if $K_X$ is not big, many important conclusions can be drawn on the birational geometry of $X$. 

A similar situation occurs when $X \subseteq \P^N$ carries a non-big Ulrich bundle, that is a bundle $\E$ with $H^i(\E(-p))=0$ for $i \ge 0, 1 \le p \le n$ and $\O_{\P(\E)}(1)$ is not big. It was shown in \cite{lo, ls} that then $X$ must be covered by linear spaces over which $\E$ is not ample. 

Moreover, at least in dimension at most $4$, one can actually classify the possible $(X, \O_X(1), \E)$, see \cite{lm, lms1}. On the other hand, if one assumes that $X$ has Picard number one, the only known varieties carrying non-big Ulrich bundles are $(\P^n, \O_{\P^n}(1))$ and quadrics $(Q_n, \O_{Q_n}(1))$. As a matter of fact, in those two cases, non-big Ulrich bundles have been classified, see \cite[Thm.~1]{lms2} or Examples \ref{sp} and \ref{qu}.

Given the above result, the absence of examples and the known results in low dimension \cite{lm, lms1}, it has been conjectured that in fact there are no non-big Ulrich bundles when the Picard number is one, except when $(X,\O_X(1)) \cong (\P^n, \O_{\P^n}(1)), (Q_n, \O_{Q_n}(1))$.

We confirm the above conjecture, as follows.

\begin{bthm} \hskip 3cm
\label{main}

Let $X  \subseteq \P^N$ be a smooth irreducible variety of dimension $n \ge 1$ such that $\rho(X)=1$. Then $X$ carries a non-big  Ulrich bundle $\E$ if and only if:
\begin{itemize}
\item [(i)] $(X,\O_X(1),\E) \cong (\P^n,\O_{\P^n}(1),\O_{\P^n}^{\oplus r})$.
\item [(ii)] $(X,\O_X(1)) \cong (Q_n,\O_{Q_n}(1))$ and $(n, \E)$ are as the following table (see Example \eqref{qu})
$$\begin{tabular}{|c|c|c|c|}
\hline
$n$ & $\E$ \\
\hline
3 & $\mathcal S$ \\
\hline
4 & $\mathcal S', \mathcal S'', \mathcal S' \oplus \mathcal S''$ \\
\hline
5 & $\mathcal S$ \\
\hline
6 & $\mathcal S', \mathcal S'', (\mathcal S')^{\oplus 2}, (\mathcal S'')^{\oplus 2}$ \\
\hline
10 & $\mathcal S', \mathcal S''$ \\
\hline
\end{tabular}.$$
\end{itemize}
\end{bthm}

\section{Notation}

Unless otherwise specified, we henceforth establish throughout the paper the following:

\begin{notation} 

\hskip 3cm

\begin{itemize} 
\item $X \subseteq \P^N$ is a smooth irreducible projective variety of dimension $n \ge 1$.
\item $\rho(X)$ is the Picard number of $X$.
\item $SX \subseteq \P^N$ is the secant variety of $X$, with dimension $s(X)$.
\item $\delta(X)=2n-s(X)+1$ is the secant defect of $X$.
\item $\nu(\L)$ is the numerical dimension of a nef line bundle $\L$ on $X$.
\item $H \in |\O_X(1)|$ is a very ample divisor. 
\item $d=H^n$ is the degree of $X$.
\item For $1 \le i \le n$, let $H_i \in |H|$ be general divisors and set $X_n:=X$ and $X_i=H_1\cap \ldots \cap H_{n-i}$, for $0 \le i \le n-1$.
\item For $x \in X$, $F_1(X,x)$ is the Fano variety of lines in $\P^N$, passing through $x$, that are contained in $X$.
\item If $\F$ is a globally generated bundle on $X$, we denote by $\xi$ the tautological line bundle of $\P(\F)=\Proj(\Sym(\F))$, $\nu(\F)= \nu(\xi)$, by $\pi : \P(\F) \to X$ the bundle map and by $\varphi=\varphi_{\xi} : \P(\F) \to \P H^0(\F)$ the associated morphism.
\item If $X$ is a Fano manifold, we set $i_X$ for its index.
\end{itemize} 
\end{notation}

For a vector bundle $\E$ on $\P^n$, we will treat the Chern classes and the Segre classes as integers. Moreover, we define $\sigma_i(\E)=s_i(\E^*)$ for the Segre classes of the dual. 

We work over the complex numbers.

\section{Two results on vector bundles on projective space}

The proof of the main theorem will proceed by studying some vector bundles associated to the fibers of the contraction morphism given by the tautological line bundle of a non-big Ulrich bundle. Since the fibers are linear spaces by \cite[Thm.~2]{ls}, we will need two technical results on vector bundles on projective spaces. 

As usual, here and in the sequel, we write Chern classes of bundles on $\P^s$ as integers. 

The first one is

\begin{lemma}
\label{primo}
Let $s\ge2$ and let $\F$ be a globally generated, non-big vector bundle on $\P^s$ with $c_1(\F)=c_1>0$.
Then:
\begin{itemize}
\item[(i)] $H^0(\F(-1))=0$.
\item[(ii)] $c_2 \ge \dfrac{c_1^2}{4}$.
\end{itemize}
\end{lemma}
\begin{proof}
Let $\pi : \P(\F) \to \P^s$ be the bundle map, let $\xi$ be its tautological line bundle and let $h=\pi^*\O_{\P^s}(1)$. Since $\F(1)$ is ample  by \cite[Cor.~6.2.12(iv)]{laz2}, we have that $\xi+h$ is ample. If $H^0(\F(-1)) \ne 0$, then $h^0(\xi-h)=h^0(\F(-1)) \ne 0$, hence there is an effective divisor $E \sim \xi-h$. Since $2\xi \sim \xi+h+E$, we get that $\xi$ is big, a contradiction. Thus, (i) is proved. To see (ii), let
$$0=\F_0 \subset \F_1 \subset \cdots \subset \F_t=\F$$
be the Harder-Narasimhan filtration, with quotients $\mathcal Q_i=\F_i/\F_{i-1}, 1 \le i \le t$ and set
$$q_i=\rk(\mathcal Q_i), \ \mu_i=\mu(\mathcal Q_i)$$
so that
$$\mu_t < \cdots < \mu_1 \ \hbox{and} \ \mu_{\rm min}(\F)=\mu_t \le \mu(\F) \le \mu_1=\mu_{\rm max}(\F).$$
Recall that the sheaves $\mathcal Q_i$ are torsion-free and slope-semistable. Fix $i \in \{1, \ldots t\}$. Since an extension of torsion-free sheaves is torsion-free, for each $j>i$, the exact sequences
$$0 \to \F_{j-1}/\F_i \to \F_j/\F_i \to \mathcal Q_j \to 0$$
show inductively that $\F/\F_i$ is torsion-free. Hence, $\F/\F_1$ is torsion-free and therefore $\F_1$ is reflexive by \cite[Prop.~1.1]{h3}, since it sits in the exact sequence 
$$0 \to \F_1 \to \F \to \F/\F_1 \to 0.$$ 
Then $q_1 \ge 2$, for otherwise $\F_1$ is a line bundle by \cite[Prop.~1.9]{h3} and, in the latter case, $\F_1 \cong \O_{\P^s}(a)$, for some $a \in \Z$. Moreover $a=\mu(\F_1)=\mu_{\rm max}(\F) \ge \mu(\F)=\frac{c_1}{\rk(\F)} >0$. Therefore we would have that $\O_{\P^s}(a) \subset \F$ and then $H^0(\F(-a)) \ne 0$, contradicting (i). Thus, $q_1 \ge 2$. Since $\F/\F_1$ is globally generated, because $\F$ is, we have that $c_1(\F/\F_1) \ge 0$ and then
$$0 < q_1 \mu(\F) \le q_1 \mu(\F_1)=c_1(\F_1)=c_1-c_1(\F/\F_1)\le c_1$$
and it follows, for $1 \le i \le t$, that
\begin{equation}
\label{mu}
\mu_i \le \mu_1=\mu(\F_1)=\frac{c_1(\F_1)}{q_1} \le \frac{c_1}{2}.
\end{equation}
The last factor $\mathcal Q_t=\F/\F_{t-1}$ is globally generated, since $\F$ is, hence $\mu_t \ge 0$. Since the factors $\mathcal Q_i$ are slope-semistable and torsion-free, Bogomolov's inequality gives, for all $i \in \{1, \ldots t\}$,
$$2q_i c_2(\mathcal Q_i)-(q_i-1)c_1(\mathcal Q_i)^2 \ge 0$$
or, equivalently, for all $i \in \{1, \ldots t\}$,
\begin{equation}
\label{ch}
\ch_2(\mathcal Q_i)=\frac{1}{2}c_1(\mathcal Q_i)^2-c_2(\mathcal Q_i) \le \frac{c_1(\mathcal Q_i)^2}{2q_i}=\frac{q_i\mu_i^2}{2}.
\end{equation}
By additivity of the Chern character we find, using \eqref{mu} and \eqref{ch},
$$\ch_2(\F)=\sum\limits_{i=1}^t \ch_2(\mathcal Q_i) \le \frac{1}{2}\sum\limits_{i=1}^t q_i\mu_i^2 \le \frac{c_1}{4} \sum\limits_{i=1}^t q_i\mu_i=\frac{c_1^2}{4}.$$
Finally $c_2=\frac{c_1^2}{2}-\ch_2(\F) \ge \frac{c_1^2}{4}$ and (ii) is proved.
\end{proof}

The second one is

\begin{lemma}
\label{bnogg}
Let $s \ge 3$ and let $\F, \F',\F''$ be vector bundles on $\P^s$, with $\F$ globally generated and non-big, $c_1(\F)=s-1$, $c_1(\F')=p \ge 2$. Suppose that we have an exact sequence
\begin{equation}
\label{ext}
 0 \to \F'(-1) \to \F \to \F'' \to 0
\end{equation}
where, either $\F''=\O_{\P^s}(p)$ or $p=2t$ and $\F''=V(t)$ for a globally generated rank $2$ bundle $V$ with $c_1(V)=t$. Then $\F'$ is not nef.
\end{lemma}
\begin{proof}
We assume by contradiction that $\F'$ is nef. Set 
$$f=\rk(\F'), \ \ell=s-1, \ u=c_2(\F) \ \hbox{and} \ v=c_3(\F).$$
By Lemma \ref{primo}(ii) and non-negativity of the Schur polynomial $c_1(\F)c_2(\F)-c_3(\F)$ \cite[Thm.~2.5]{dps}, we have
\begin{equation}
\label{uv}
u \ge \frac{\ell^2}{4}, \ v \le \ell u.
\end{equation}
Note that $H^0(\F''(-1)) \ne 0$, because either $\F''(-1)=\O_{\P^s}(p-1)$ or  $\F''(-1)=V(t-1)$ and $t=\frac{p}{2} \ge 1$. Since $H^0(\F(-1))=0$, \eqref{ext} gives that $H^1(\F'(-2)) \ne 0$. We observe that this implies that
\begin{equation}
\label{lp}
s-1=\ell \le p.
\end{equation}
Indeed, let $\xi$ be the tautological bundle on $\P(\F')$ and let $h=\pi^*H$. Since $\F'(1)$ is ample by \cite[Cor.~6.2.12(iv)]{laz2}, it follows that $\xi+h$ is ample. If $\ell \ge p+1$, then 
$$\xi-2h=K_{\P(\F')}+(f+1)\xi+(\ell-p)h$$
with $(f+1)\xi+(\ell-p)h=\xi+h+f\xi+(\ell-p-1)h$ is ample, because $\xi$ and $h$ are nef. 
Then, Kodaira vanishing implies that $0=H^1(\xi-2h)=H^1(\F'(-2))$, a contradiction. Hence \eqref{lp} is proved.

Next, we compute Chern classes of $\F'$. To ease calculations we set
$$T(x)=\frac{x^3}{3}+\frac{p+2}{4}x^2+ \frac{3p+4}{6}x$$
and we observe that, being $T'(x) \ge 0$ for $x \ge 0$, it follows that $T(x)$ is an increasing function.

Assume now that we are in the case $\F''=\O_{\P^s}(p)$.

Computing first Chern classes in \eqref{ext} gives $f=2p-\ell$ and it follows that
$$c_1(\F'(-1))=\ell-p, \ c_2(\F'(-1))=u-p\ell+p^2, \ c_3(\F'(-1))=v-pu+p^2\ell-p^3.$$
By the above, we get, using $f=2p-\ell$,
$$c_2(\F')=u+p^2-\frac{\ell(\ell+1)}{2}$$
and
$$c_3(\F')=v+(p-\ell-2)u+\frac{p^3-3p^2-p}{3}-p\binom{\ell+1}{2}+2\binom{\ell+2}{3}.$$
It follows that
\begin{equation}
\label{sig1}
\sigma_3(\F')=v-(l+p+2)u-\frac{p(p+1)(2p+1)}{3}+p\binom{\ell+1}{2}+2\binom{\ell+2}{3}.
\end{equation}
Using $v \le lu$ and then $u\ge l^2/4$, we obtain
$$\begin{aligned}
\sigma_3(\F') & \le -(p+2)u+p\binom{\ell+1}{2}+2\binom{\ell+2}{3}-\frac{p(p+1)(2p+1)}{3} \\
& \le -\frac{(p+2)l^2}{4}+p\binom{\ell+1}{2}+2\binom{\ell+2}{3}-\frac{p(p+1)(2p+1)}{3} \\
&=\frac{l^3}{3}+\frac{(p+2)l^2}{4}+\frac{(3p+4)l}{6}-\frac{p(p+1)(2p+1)}{3} \\
&=T(l)-\frac{p(p+1)(2p+1)}{3} \le T(p)-\frac{p(p+1)(2p+1)}{3}=\frac{p(4-p^2)}{12}.
\end{aligned}$$
If $s \ge 4$, then $p \ge 3$ by \eqref{lp}. Choose a linear space $P \subset \P^s$ with $\dim P=3$. Then 
$$\sigma_3(\F'_{|P}) \le \frac{p(4-p^2)}{12}<0$$ 
a contradiction, by \cite[Rmk.~2.2]{lm}, since $\F'_{|P}$ is nef. If $s=3$ we have that $\ell=2$. Moreover, being $\F$ non big, we have that $0=\sigma_3(\F)=c_1(\F)^3-2c_1(\F)c_2(\F)+c_3(\F)=8-4u+v$, so that $v=4u-8 \ge 0$, hence $u \ge 2$. Replacing into \eqref{sig1}, we find that
$$\sigma_3(\F')=\frac{8p}{3}-p^2-\frac{2p^3}{3}-pu \le -\frac{p(p+2)(2p-1)}{3}<0$$
again a contradiction as above. Thus the lemma is proved in the case $\F''=\O_{\P^s}(p)$.

Finally, assume now that we are in the case $\F''=V(t)$. We set 
$$w=c_2(V)$$ 
so that, $c_1(V(t))=3t, c_2(V(t))=w+2t^2$. Comparison of first Chern classes in \eqref{ext}
gives
$$f=5t-\ell.$$
We therefore obtain
$$c_1(\F'(-1))=\ell-3t, \ c_2(\F'(-1))=u-3t \ell+7t^2-w, \ c_3(\F'(-1))=v-3tu+\ell(7t^2-w)+6tw-15t^3.$$
By the above, we get 
$$c_1(\F')=2t, \ c_2(\F')=u-w+\frac{9t^2+t-\ell(\ell+1)}{2}$$
and
$$c_3(\F')=\frac{2\ell}{3}+\ell^2+\frac{\ell^3}{3}-\frac{4t}{3}-\ell t-\ell^2 t-4t^2+\frac{10t^3}{3}-2u-\ell u+2tu+v+2w+tw.$$
It follows that
\begin{equation}
\label{sig3}
\sigma_3(\F')=v-(\ell+2t+2)u+(5t+2)w-\frac{2t-\ell}{3}\left(\ell^2+5\ell t+3\ell+10t^2+9t+2\right).
\end{equation}
Using \eqref{uv} and \eqref{lp}, we get
\begin{equation}
\label{sig2}
\begin{aligned}
\sigma_3(\F') & \le -(2t+2)u+(5t+2)w-\frac{2t-\ell}{3}\left(\ell^2+5\ell t+3\ell+10t^2+9t+2\right) \\
& \le -\frac{t+1}{2}\ell^2+(5t+2)w-\frac{2t-\ell}{3}\left(\ell^2+5\ell t+3\ell+10t^2+9t+2\right) \\
& = T(\ell)+(5t+2)w-\frac{2t(5t+2)(2t+1)}{3} \\
& \le T(2t)+(5t+2)w-\frac{2t(5t+2)(2t+1)}{3} \\
& = (5t+2)w-2t^2(t+1).
\end{aligned}
\end{equation}
Now, global generation of $V$ gives, using the non-negativity of the Schur polynomials $c_1(V)^3-2c_1(V)c_2(V)$ and, for $s \ge 4$, $c_1(V)^4-3c_1(V)^2c_2(V)+c_2(V)^2$, the conditions
$$t^3-2tw \ge 0 \ \hbox{and, for} \ s \ge 4, \ t^4-3t^2w+w^2 \ge 0.$$
If $s \ge 4$, we get $w \le (3-\sqrt5)t^2/2 < 2t^2/5$ and replacing into \eqref{sig2} we find $\sigma_3(\F') \le -6t^2/5<0$, a contradiction. If $s=3$, we have, as above, that $\ell=2, v = 4u-8$ and $u \ge 2$. Replacing for $\ell$ and $v$ into \eqref{sig3}, and using $u \ge 2$ and $w \le \frac{t^2}{2}$, we find that
$$\sigma_3(\F') \le -\frac{2t(2t-1)(5t+7)}{3}-2tu+(5t+2)w \le -\frac{t(25t^2+30t-4)}{6}<0$$
again a contradiction. This completes the proof of the lemma.
\end{proof}

\section{Two results on subcanonical varieties}

In this section we prove two useful results on subcanonical varieties.

First, we need the following, probably well-known, fact. 

\begin{lemma}
\label{nf}
Let $Y$ be a smooth irreducible variety and let $f : Y \to T$ be a projective surjective morphism, with $T$ an integral variety such that $0< \dim T < \dim Y$. Let $F$ be a smooth scheme-theoretic fiber of $f$. Then $N_{F/Y}^*$ is globally generated and non-big.
\end{lemma}
\begin{proof}
Let $F$ be the fiber over $t \in T$ and let $g : F \to \{t\}$. It follows by \cite[Lemma 29.32.4, \href{https://stacks.math.columbia.edu/tag/0473}{Tag 0473}]{st} that there is a surjection $g^* \mathcal C_{\{t\}/T} \to \mathcal C_{F/Y}=N_{F/Y}^*$. On the other hand, $\mathcal C_{\{t\}/T} = \I_{\{t\}/T}/\I_{\{t\}/T}^2 \cong \O_{\{t\}}^{\oplus a}$, where $a= \dim_{\C(t)} m_t/m_t^2$, and therefore we get a surjection $\O_F^{\oplus a} \cong g^*\O_{\{t\}}^{\oplus a} \to N_{F/Y}^*$, showing that $N_{F/Y}^*$ is globally generated. Consider the blow up $\widetilde T$ of $T$ in $t$, with exceptional divisor $E_1$ and the blow up $p : \widetilde Y \to Y$ of $Y$ along $F$, with exceptional divisor $E$. Since 
$$\I_{F/Y}=\Im \{f^* \I_{\{t\}/T} \to \O_Y\}$$ 
and 
$$\O_{\widetilde Y}(-E)=\I_{F/Y} \cdot \O_{\widetilde Y}=\Im \{p^* \I_{F/Y} \to \O_{\widetilde Y}\}=\Im \{(f \circ p)^*\I_{\{t\}/T} \to \O_{\widetilde Y}\}=\I_{\{t\}/T} \cdot \O_{\widetilde Y}$$
it follows by the universal property of blowing up that there is a morphism $\widetilde f : \widetilde Y \to \widetilde T$ over $T$, giving rise to a commutative diagram
$$\begin{tikzcd}
E \arrow[r, hook, "j"] \arrow[d, "h"'] & \widetilde Y \arrow[d, "\widetilde f"] \\ E_1 \arrow[r, hook, "i"] & \widetilde T.
\end{tikzcd}$$
Since $\P(N_{F/Y}^*) \cong E$, we have that 
$$\O_{\P(N_{F/Y}^*)}(1) \cong \O_E(-E)=j^* \O_{\widetilde Y}(-E) \cong j^*{\widetilde f}^* \O_{\widetilde T}(-E_1) \cong h^*(i^* \O_{\widetilde T}(-E_1))$$
hence $\O_{\P(N_{F/Y}^*)}(1)$ is not big, since $\dim E_1 = \dim \widetilde T - 1 < \dim \widetilde Y - 1 = \dim E$.
\end{proof}

The first result implies that $b_4(X) = 1$ on some Fano manifolds.

\begin{lemma}
\label{b4}
Let $X \subset \P^N$ be a smooth irreducible variety of dimension $n \ge 2$ such that $-K_X=eH$, with either $e > \max\{\frac{2n}{3},\frac{n+4}{2}\}$ or with secant defect $\delta(X) \ge 5$. Then $H^4(X;\Z) \cong \Z H^2$.
\end{lemma}
\begin{proof}
We can assume that $X$ is nondegenerate. If $X$ is a linear space, the conclusion is immediate; assume otherwise. Set $s=s(X)$ and $\delta=\delta(X)=2n+1-s$. We first show that, if $e > \max\{\frac{2n}{3},\frac{n+4}{2}\}$, then $\delta \ge 5$. Indeed, we have $i_X \ge e > \frac{2n}{3}$, hence it follows from \cite[Cor.~2]{ir2} that 
$X$ is connected by degenerate conics.
 By \cite[Lemma 5(i)]{ir2}, we have that $\delta \ge 2$ and it follows by \cite[Cor.~3(i)]{ir2} that $\delta \ge 2i_X-n \ge 2e-n> 4$, so that $\delta \ge 5$. Assume now that $\delta \ge 5$. Then
$$2n-s=\delta-1 \ge 4.$$
If $s<N$, take a general linear space $\Lambda \cong \P^{N-s-1}$ disjoint from $SX$. Projection from $\Lambda$ gives a closed embedding $X \subset \P^s$. If $s=N$, use the original embedding. Then, the integral Barth-Larsen theorem (see for example \cite[Thm.~2.2(b)]{h2}) implies that
$$H^4(\P^s;\Z) \mapright{\cong} H^4(X;\Z)$$
and therefore $H^4(X;\Z) \cong \Z H^2$. 
\end{proof}

The second result, is about planes contained in a subcanonical variety.

\begin{lemma}
\label{plane}
Let $X \subset \P^N$ be a smooth irreducible $n$-dimensional variety of codimension $c \ge 1$ with $-K_X=eH$, for some $e \in \Z, H \in |\O_X(1)|$. Let $P \subset X \subset \P^N$ be a linear plane. Then, $N_{P/X}^*(1)$ is globally generated,
\begin{equation}
\label{pl1}
c_1(N_{P/X}^*(1))=n+1-e, \ \hbox{and, if} \ c \ge 2, \ \sigma_2(N_{P/X}^*(1)) \ge n-e.
\end{equation}
\end{lemma}
\begin{proof}
We can assume that $X$ is nondegenerate. For simplicity of notation we write $N_P=N_{P/\P^N}, N_X=N_{X/\P^N}$ and $\I_P=\I_{P/\P^N}$. From the exact sequence 
$$0 \to T_P \to {T_X}_{|P} \to N_{P/X} \to 0$$
we get that $c_1(N_{P/X})=c_1(X)_{|P}-3=e-3$, hence
$$c_1(N_{P/X}^*(1))=n+1-e.$$
Assume now that $c \ge 2$. Since 
$$N_P^*(1) \cong H^0(\I_P(1)) \otimes \O_P$$
we have an exact sequence
\begin{equation}
\label{n1}
0 \to N_X^*(1)_{|P} \to H^0(\I_P(1)) \otimes \O_P \to N_{P/X}^*(1) \to 0
\end{equation}
showing that $c_1(N_X^*(1)_{|P})=-c_1(N_{P/X}^*(1))$ and
\begin{equation}
\label{n2}
c_2(N_X^*(1)_{|P})=-c_2(N_{P/X}^*(1))-c_1(N_X^*(1)_{|P})c_1(N_{P/X}^*(1))=c_1(N_{P/X}^*(1))^2-c_2(N_{P/X}^*(1))=\sigma_2(N_{P/X}^*(1)).
\end{equation}
Note that $c_2(N_X(-1)_{|P}) \ge 0$ since, as is well known, $N_X(-1)$ is a globally generated rank $c$ bundle, hence so is $N_X(-1)_{|P}$. We now prove that 
\begin{equation}
\label{c2} 
c_2(N_X(-1)_{|P}) > 0.
\end{equation}
Suppose that $c_2(N_X(-1)_{|P})=0$. Choosing a general subspace $V \subset H^0(N_X(-1)_{|P})$ with $\dim V =c-1$, we get an exact sequence
\begin{equation}
\label{sp1}
0 \to \O_P^{\oplus (c-1)} \mapright{\varphi} N_X(-1)_{|P} \to \I_{Z/P}(n+1-e) \to 0
\end{equation}
where $Z=D_{c-2}(\varphi)$ is the degeneracy locus. It follows by \cite[Statement (folklore)(i), \S 4.1]{ba} that either $Z=\emptyset$ or $Z$ is of pure codimension $2$. But, in the latter case, $[Z]=[D_{c-2}(\varphi)]=c_2(N_X(-1)_{|P})=0$, a contradiction. Therefore  $Z=\emptyset$ and \eqref{sp1} becomes
$$0 \to \O_P^{\oplus (c-1)} \to N_X(-1)_{|P} \to \O_P(n+1-e) \to 0.$$
Since $\Ext^1(\O_P(n+1-e),\O_P^{\oplus (c-1)}) \cong H^1(\O_P(-n-1+e))^{\oplus (c-1)}=0$, the above sequence splits and
$$N_X(-1)_{|P} \cong \O_P^{\oplus (c-1)} \oplus \O_P(n+1-e).$$
We claim that the direct summand $\O_P^{\oplus (c-1)}$ gives rise to a linear space $\Lambda \subset \P^N$ of dimension $n+1$ such that, denoting by $\mathbb T_xX$ the embedded tangent space to $X$ at $x$,
\begin{equation}
\label{z1}
\mathbb T_xX \subset \Lambda, \ \hbox{for every} \ x \in P.
\end{equation}
Indeed, after dualizing the above splitting, the direct summand $\O_P^{\oplus (c-1)}$ determines a vector subspace $W \subseteq H^0({N_X^*(1)}_{|P}) \subseteq H^0(\I_P(1))$ (by \eqref{n1}), with $\dim W = c-1$ and $W \cong W \otimes \C(x) \subset (N_X^*(1)_{|P})(x)$ for every $x \in P$. As is well known, $(N_X^*(1)_{|P})(x)$ is, via the Euler sequence, the space of hyperplanes containing $\mathbb T_xX$. It follows that, for every $x \in P$,
$$\mathbb T_xX \subseteq \Lambda:=\bigcap\limits_{H \in W} H.$$
Now $\dim \Lambda = N-\dim W=N-c+1=n+1$ and \eqref{z1} is proved. But then, 
$$P \subseteq \{x \in X : \mathbb T_xX \subseteq \Lambda\}$$
and Zak's theorem on tangencies \cite[Thm.~3.4.17]{la} implies the contradiction $2= \dim P \le \dim \Lambda -n=1$. This proves \eqref{c2}. If $c \ge 3$, since $N_X(-1)_{|P}$ is globally generated of rank $c$, we can choose a general subspace $V' \subset H^0(N_X(-1)_{|P})$ with $\dim V' =c-2$, giving rise to an exact sequence
$$0 \to \O_P^{\oplus (c-2)} \mapright{\psi} N_X(-1)_{|P} \to \F \to 0$$
where $\F$ is the quotient. As above, since $c_3(N_X(-1)_{|P})=0$, it follows by \cite[Statement (folklore)(i), \S 4.1]{ba} that $D_{c-3}(\psi)=\emptyset$, and therefore $\F$ is a globally generated rank $2$ bundle on $P$ with the same Chern classes of $N_X(-1)_{|P}$. 
If $c=2$, we set $\F=N_X(-1)_{|P}$. It follows by \eqref{c2} and \cite[Lemma 2.4]{e}, that 
$$c_2(\F) \ge c_1(\F)-1.$$ 
Now,
$$c_1(\F)=c_1(N_X(-1)_{|P})=-c_1(N_X^*(1)_{|P})=c_1(N_{P/X}^*(1))=n+1-e$$
and, by \eqref{n2},
$$c_2(\F)=c_2(N_X(-1)_{|P})=c_2(N_X^*(1)_{|P})=\sigma_2(N_{P/X}^*(1))$$
and it follows that 
$$\sigma_2(N_{P/X}^*(1))=c_2(\F) \ge c_1(\F)-1=n-e$$
so that \eqref{pl1} is proved.
\end{proof} 

\section{Numerical dimension, Segre classes and a fact about Ulrich bundles}

We give a simple result about numerical dimension that will allow inductive proofs. We include a proof for completeness' sake.

\begin{lemma}
\label{sez}
Let $X \subset \P^N$ be a smooth irreducible variety of dimension $n$ and let $\E$ be a rank $r$ globally generated bundle on $X$. Then 
\begin{itemize}
\item [(i)] $\nu(\E)=\max\{k \ge r-1: s_{k-r+1}(\E) \ne 0\}$.
\item [(ii)] If $\E$ is not big and $\nu=\nu(\E)$ (so that $r-1 \le \nu \le n+r-2$), then $\E_{|X_{\nu-r+1}}$ is big and $\E_{|X_i}$ is not big for $\nu-r+2 \le i \le n$.
\end{itemize}
\end{lemma}
\begin{proof}
Recall that, by definition of Segre classes (see for example \cite[Ch.~10]{eh}, where they use $\Proj(\Sym(\E^*))$), $s_i(\E^*)=(-1)^i s_i(\E)=\pi_*(\xi^{r-1+i})$ for $0 \le i \le n$. Now $\xi^l=0$ for $l \ge \nu+1$, hence $s_{l-r+1}(\E^*)=\pi_*(\xi^{r-1+(l-r+1)})=0$ for $l \ge \nu+1$. On the other hand, $\xi^{\nu} \ne 0$, hence there is $\alpha \in A^{n+r-1-\nu}(\P(\E))$ such that $\xi^{\nu} \alpha \ne 0$ (for example $\alpha$ a complete intersection of very ample divisors). Writing $\alpha = \sum_{i=0}^{r-1} \xi^i \pi^* \alpha_i$ we obtain, as intersection numbers, 
$$0 \ne \xi^{\nu} \alpha = \xi^{\nu} \pi^* \alpha_0  = \pi_*(\xi^{\nu}) \alpha_0 =s_{\nu-r+1}(\E^*) \alpha_0$$
so that $s_{\nu-r+1}(\E) \ne 0$ and (i) is proved. To see (ii), note that $\pi^* H^{n-i} = \P(\E_{|X_i})$, hence, if $i \ge \nu-r+2$ we have that 
$$(\xi_{\E_{|X_i}})^{i+r-1} = \xi^{i+r-1}\pi^* H^{n-i}=0$$
hence $\E_{|X_i}$ is not big. On the other hand, 
$$(\xi_{\E_{|X_{\nu-r+1}}})^{\nu} = \xi^{\nu}\pi^* H^{n+r-1-\nu} \ne 0$$
for otherwise $\pi_*(\xi^{\nu})H^{n+r-1-\nu}=0$, hence $\pi_*(\xi^{\nu})=0$, since $\pi_*(\xi^\nu)$ is represented by an effective cycle and $H$ is ample, so that $s_{\nu-r+1}(\E)=0$, contradicting (i). This proves (ii).
\end{proof}

In the case of Ulrich bundles, we set up the necessary notation and establish a useful fact.

\begin{defi-not}
\label{basica}
Let $X \subset \P^N$ be a smooth irreducible variety and let $\E$ be an Ulrich bundle on $X$. We denote by 
$$\pi: Y=\P(\E) \to X$$
the bundle map, with tautological line bundle $\xi$ and, for $H \in |\O_X(1)|$, we set $h=\pi^*H$. We denote by
$$\varphi=\varphi_{\xi} : Y \to Z:=\varphi(Y) \subseteq \P H^0(\E)$$ 
the map associated to $\xi$.
\end{defi-not}

\begin{remark}
\label{basic2}
Note that, since $\E$ is Ulrich, it is $0$-regular, hence globally generated. Therefore so is $\xi$ and this defines $\varphi$. We will often use the fact, proved in \cite[Thm.~2]{ls}, that, for any scheme-theoretic fiber $F$ of $\varphi$, we have that $F \cong \P^m$ for some $m \in \Z$, $\xi_{|F} \cong \O_F, h_{|F} \cong \O_{\P^m}(1)$ and $\pi_{|F} : F \to \pi(F) \subset X$ is an isomorphism, making $\pi(F)$ a linear space contained in $X \subset \P^N$. In the sequel, we will identify $F$ with $\pi(F)$ and we will write accordingly $F \subset Y$ or $F \subset X$.
\end{remark}

We have 
\begin{lemma}
\label{ul}
Let $X \subset \P^N$ be a smooth irreducible variety of dimension $n \ge 2$ and let $\E$ be an Ulrich bundle on $X$. We have:
\begin{itemize}
\item [(i)] $H^0(T_X(-2))=0$.
\item [(ii)] $H^0(T_X \otimes \E(-2))=0$.
\item [(iii)] If there is a linear plane $P \subset X \subset \P^N$, then $H^1(\E_{|P}(-1))=H^2(\E_{|P}(-1))=0$.
\end{itemize}
\end{lemma}
\begin{proof}
The twisted Euler sequence
\begin{equation}
\label{e}
0 \to \O_X(-2) \to \O_X(-1)^{\oplus(N+1)} \to {T_{\P^N}}_{|X}(-2) \to 0
\end{equation}
gives, by Kodaira vanishing, that $H^0({T_{\P^N}}_{|X}(-2))=0$ and therefore, since $T_X \hookrightarrow {T_{\P^N}}_{|X}$, (i) follows. Tensoring \eqref{e} with $\E$ gives an exact sequence
$$0 \to \E(-2) \to \E(-1)^{\oplus(N+1)} \to {T_{\P^N}}_{|X} \otimes \E(-2) \to 0$$
and since $\E$ is Ulrich, we have $H^0(\E(-1))=H^1(\E(-2))=0$, so that $H^0({T_{\P^N}}_{|X} \otimes \E(-2))=0$ and therefore (ii) again follows. To see (iii), consider the strand of the resolution of $\E$ as a $\O_{\P^N}$-module \cite[Prop.~2.1]{es}, \cite[Thm.~2.3(i)]{b}),
$$\ldots \to \O_{\P^N}(-1)^{\oplus \beta_1} \to \O_{\P^N}^{\oplus \beta_0} \to \E \to 0.$$
Tensoring by $\O_P(-1)$, we get an exact sequence 
$$\ldots \to \O_P(-2)^{\oplus \beta_1} \to \O_P(-1)^{\oplus \beta_0} \to \E_{|P}(-1) \to 0$$
that we can split into two exact sequences
\begin{equation}
\label{k1}
0 \to \K \to \O_P(-1)^{\oplus \beta_0} \to \E_{|P}(-1) \to 0
\end{equation}
and 
$$\O_P(-2)^{\oplus \beta_1} \to \K \to 0.$$
Now, the latter shows that $\K(2)$ is globally generated and we claim that 
\begin{equation}
\label{k}
H^0(\K^*(-3))=0.
\end{equation}
Indeed, since $H^0(\K^*(-3)) \subseteq H^0(\K^*(-2))$, \eqref{k} certainly holds if $H^0(\K^*(-2))=0$. If $H^0(\K^*(-2)) \ne 0$ we have that $s:=h^0((\K(2))^*) \ge 1$, hence, as is well known, $\K(2) \cong \O_P^{\oplus s} \oplus \G$, where $\G$ is a vector bundle on $P$ with $H^0(\G^*)=0$. But then $\K^*(-3) \cong \O_P(-1)^{\oplus s} \oplus \G^*(-1)$ and therefore $H^0(\K^*(-3))=0$, proving \eqref{k}.
It follows from \eqref{k1} that $H^1(\E_{|P}(-1))=0$, since $H^1(\O_P(-1)^{\oplus \beta_0})=0$ and, by Serre's duality, $H^2(\K) \cong H^0(\K^*(-3))^*=0$ by \eqref{k}. Moreover, still \eqref{k1}, and $H^2(\O_P(-1)^{\oplus \beta_0})=H^3(\K)=0$, imply that $H^2(\E_{|P}(-1))=0$ and (iii) is proved.
\end{proof} 

\section{Known examples and low dimension}

The following are known examples of non-big Ulrich bundles on a variety $X$ with $\rho(X)=1$.

\begin{example}
\label{sp}
On $(\P^n, \O_{\P^n}(1))$ we know (see for example \cite[Prop.~2.1]{es}, \cite[Thm.~2.3]{b}) that the only rank $r$ Ulrich bundle is $\O_{\P^n}^{\oplus r}$ and it is not big (for example by \cite[Rmk.~2.2]{lm}, because $c_1(\E)^n=0$). 
\end{example}

\begin{example}
\label{qu}
On $(Q_n, \O_{Q_n}(1)), n \ge 2$, we know by \cite[Thm.~1]{lms2} that the only non-big Ulrich bundles are the ones in the following table (where $\mathcal S, \mathcal S', \mathcal S''$ are the spinor bundles, see \cite[Def.~3.1]{lms2}):
$$\begin{tabular}{|c|c|c|c|}
\hline
$n$ & $\E$ \\
\hline
2 & $(\mathcal S')^{\oplus r}, (\mathcal S'')^{\oplus r}, r \ge1$ \\
\hline
3 & $\mathcal S$ \\
\hline
4 & $\mathcal S', \mathcal S'', \mathcal S' \oplus \mathcal S''$ \\
\hline
5 & $\mathcal S$ \\
\hline
6 & $\mathcal S', \mathcal S'', (\mathcal S')^{\oplus 2}, (\mathcal S'')^{\oplus 2}$ \\
\hline
10 & $\mathcal S', \mathcal S''$ \\
\hline
\end{tabular}.$$ 
\end{example}

In fact when $X \cong \P^n$ or $X \cong Q_n$ with $n \ge 2$, we don't need to worry about the polarization for the non-big Ulrich bundle.

\begin{lemma}
\label{pq}
Let $X \subset \P^N$ be a smooth irreducible variety of dimension $n \ge 1$ and let $\E$ be a rank $r$ non-big Ulrich bundle on $X$.
\begin{itemize}
\item [(i)] If $X \cong \P^n$, then $(\O_X(1), \E) \cong (\O_{\P^n}(1), \O_{\P^n}^{\oplus r})$.
\item [(ii)] If $X \cong Q_n$ and $n \ge 2$, then $\O_X(1) \cong \O_{Q_n}(1)$ and $(n, \E)$ are as in Example \ref{qu}.
\end{itemize}
\end{lemma}
\begin{proof}
We have that $\O_X(1) \cong \O_{\P^n}(t)$ or $\O_{Q_n}(t)$ for some $t \ge 1$. On the other hand, $\E$ being non-big implies, by \cite[Thm.~1]{lo}, that $t=1$ and we conclude by Examples \ref{sp} and \ref{qu}.
\end{proof}

In the case of a curve, we have

\begin{lemma}
\label{curva}
Let $X \subseteq \P^N$ be a smooth irreducible curve and let $\E$ be an Ulrich bundle on $X$. Then $\E$ is not big if and only if $(X, \O_X(1), \E)=(\P^1, \O_{\P^1}(1), \O_{\P^1}^{\oplus r})$.
\end{lemma}
\begin{proof}
If $(X, \O_X(1), \E)=(\P^1, \O_{\P^1}(1), \O_{\P^1}^{\oplus r})$, we know that $\E$ is not big by Example \ref{sp}. Vice versa, assume that $\E$ is not big. If $g$ is the genus of $X$, then $\deg \E=r(d+g-1)$ by \cite[Prop.~2.3]{ch}. Now, it follows from \cite[Rmk.~2.2]{lm} that $0=s_1(\E^*)=c_1(\E)=r(d+g-1)$, so that $d=1$ and $g=0$, hence $(X, \O_X(1))=(\P^1, \O_{\P^1}(1))$ and $\E = \O_{\P^1}^{\oplus r}$ by Example \ref{sp}. 
\end{proof}

In dimension at most $4$ we can use known results.

\begin{thm}
\label{nle4}
Let $X \subset \P^N$ be a smooth irreducible variety of dimension $n \le 4$ with $\rho(X)=1$ and let $\E$ be an Ulrich bundle on $X$. Then $\E$ is not big if and only if $(X, \O_X(1),\E)$ is as in Examples \ref{sp} and \ref{qu} (with $n \le 4$).
\end{thm}
\begin{proof}
Follows from Lemma \ref{curva} and \cite[Cor.~1(i)]{lms1}.
\end{proof}

\section{Varieties with $\rho(X)=1$ and a non-big Ulrich bundle}

We start by collecting some important consequences on the geometry of $X$ in the presence of a non-big Ulrich bundle.

\begin{lemma}
\label{nongg}
Let $X \subset \P^N$ be a smooth irreducible variety of dimension $n \ge 1$ and let $\E$ be a rank $r$ non-big Ulrich bundle on $X$. Then $K_X+\det \E$ is not globally generated and if $F$ is a general fiber of $\varphi=\varphi_{\xi}$, then
$F \cong \P^k$, with $k=n+r-1-\nu(\E) \ge 1$ and
$$(\pi^*(K_X+\det \E))_{|F} = \O_{\P^k}(-k-1).$$
\end{lemma}
\begin{proof}
The fact that $F \cong \P^k$, with $k=n+r-1-\nu(\E)$ follows from Remark \ref{basic2} and $k \ge 1$ since $\E$ is not big. Now,
$$\O_{\P^k}(-k-1)=K_{\P^k}=K_F={K_Y}_{|F}=(-r\xi+\pi^*(K_X+\det \E))_{|F}=(\pi^*(K_X+\det \E))_{|F}$$
and therefore $K_X+\det \E$ is not globally generated.
\end{proof}

The following result will show already strong constraints on $X$.
\begin{prop} 
\label{uno}
Let $X  \subseteq \P^N$ be a smooth irreducible variety of dimension $n \ge 2$, degree $d \ge 2$ and such that $\rho(X)=1$. Let $\E$ be a rank $r$ non-big Ulrich bundle on $X$. Then the following hold:
\begin{itemize}
\item [(i)] $X$ is Fano and is covered by lines.
\item [(ii)] $\Pic(X) \cong \Z \O_X(1)$.
\item [(iii)] $-K_X=i_XH$ with $i_X \ge n+r+1-\nu(\E) \ge 3$.
\item [(iv)] $c_1(\E) = \frac{r}{2}(n+1-i_X)H$ and $c_1(\E)^n>0$. 
\item [(v)] $r \ge 2$.
\item [(vi)] $\nu(\E) \ge \frac{n}{2}+r-1$.
\item [(vii)] $\Nef(\P(\E))=\Pseff(\P(\E))=\{a\xi+bh, a \ge 0, b \ge 0 \}$.
\item [(viii)] $\xi+h$ is very ample.
\item [(ix)] $\E$ is a Fano bundle (that is $\P(\E)$ is Fano) and 
$$\nu(\E)=\frac{(r+2)(n-i_X)+3r}{2}$$
$$\tau(\E)=\frac{r}{2}(n+1-i_X)$$
where $\tau(\E)$ is the only real number such that $-K_{\P(\E)} + \pi^*(K_X) + \tau(\E)h$ is nef non ample.
\item [(x)] $$\frac{2i_X-n-2}{n+1-i_X} \le r \le \frac{2i_X-4}{n+1-i_X}.$$ 
Equivalently, 
\begin{equation}
\label{indi}
\frac{r(n+1)+4}{r+2} \le i_X \le \frac{r(n+1)+n+2}{r+2}.
\end{equation}
Moreover, if $X \not\cong Q_n$, then $r \le n-3$ and if equality holds, $i_X=n-1$, that is $(X, \O_X(1))$ is a Del Pezzo manifold.
\end{itemize}
\end{prop}
\begin{proof}
Since $d \ge 2$ we have that $(X, \O_X(1)) \ne (\P^n, \O_{\P^n}(1))$. By hypothesis there is an ample divisor $A$ such that $H \equiv hA$ and $K_X \equiv aA$ for some integers $h, a$ with $h \ge 1$. Note that through a general point $x \in X$ there is a line $L \subset X$: This follows from \cite[Thm.~1]{lo} or from \cite[Thm.~2]{ls}. Hence
$$1 = H \cdot L = h A \cdot L$$ so that $h=1$. As is well-known (see for example \cite[Prop.~2.3.9]{ru}), if $x \in X$ is a general point, we have that $0 \le \dim_{[L]} F_1(X,x) = -K_X \cdot L-2 = -a - 2$, hence $a <0$ and $X$ is Fano. This proves (i). Consequently, the numerical and linear equivalence for divisors coincide on $X$ and then $H \sim A$, that is $H$ generates $\Pic(X)$, so that (ii) is proved. 

Now $\det \E=\frac{r}{2}(K_X+(n+1)H)$ by \cite[Lemma 3.2]{lo} and (iv) follows since, as is well known, $i_X \le n$ unless $X \cong \P^n$, but then $d=1$ by (ii). By (iv) and \cite[Cor.~2]{ls}, we have that $\dim F_1(X,x) \ge n+r-1-\nu(\E)$. Since $\dim F_1(X,x)=-K_X \cdot L-2=i_X-2$, we get that $i_X \ge n+r+1-\nu(\E)$. On the other hand, $\E$ is not big, hence $\nu(\E) \le \dim \P(\E)-1=n+r-2$ and this gives (iii). 

To see (v), if $r=1$, we have by (ii) that $\E=\O_X(l)$ is Ulrich, and this implies that $(X, \O_X(1))=(\P^n, \O_{\P^n}(1))$ (see for example \cite[Lemma 4.2(vi)]{aclr}), a contradiction. Next, (vi) follows from \cite[Thm.~2]{lms2}. To see (vii), by (v) we can choose a line $L_{\pi} \subseteq \P(\E(x)) \cong \P^{r-1}$, so that $\xi \cdot L_{\pi}=1, h \cdot L_{\pi}=0$. On the other hand, consider $\varphi=\varphi_{\xi} : Y \to Z$. It follows from Remark \ref{basic2} that for any fiber $F$ of $\varphi$, if $\Gamma \subseteq \pi(F)$ is a line and $L_{\varphi} \subseteq F \subset \P(\E)$ is the corresponding curve under the isomorphism $F \cong \pi(F))$, then $\xi \cdot L_{\varphi} =0, h \cdot L_{\varphi} =1$. Now, since $\xi$ and $h$ are nef, we get that $a\xi+bh$ is nef when $a \ge 0, b \ge 0$. On the other hand, if $a\xi+bh$ is nef, we get that $a \ge 0, b \ge 0$ by intersecting with $L_{\pi}$ and $L_{\varphi}$. Next, we prove that $a\xi+bh$ is big if and only if $a>0, b>0$. Clearly this concludes the proof of (vii). If $a>0, b>0$, then $a\xi+bh$ is nef and 
$$(a\xi+bh)^{n+r-1}=\sum\limits_{j=0}^{n+r-1} \binom{n+r-1}{j}a^jb^{n+r-1-j}\xi^j h^{n+r-1-j}.$$
In the above sum, all terms are non-negative, since $\xi$ and $h$ are nef. On the other hand, $\xi^{r-1} h^n=H^n>0$ and therefore $(a\xi+bh)^{n+r-1}>0$, hence $a\xi+bh$ is big. Finally, assume that $a\xi+bh$ is big. We have that $a\xi+bh \equiv A+E$ with $A$ ample and $E$ effective. Choosing a line $L_{\pi}$ as above and such that $L_{\pi} \not\subset \Supp(E)$, we get that $E \cdot L_{\pi} \ge 0$, hence $a = (a\xi+bh) \cdot L_{\pi} =A\cdot L_{\pi}+E\cdot L_{\pi} > 0$. Now, if $b \le 0$, we have that also $a\xi=a\xi+bh-bh$ is big, a contradiction since $a>0$ and $\xi$ is not big. This proves (vii). 

As for (viii), note that $\xi+h$ is the tautological line bundle of $\E(H)$, which is strongly very ample because $\E$, being an Ulrich bundle, is globally generated. Therefore $\E(H)$ is very ample, that is $\xi+h$ is very ample, by \cite[Prop.~4.1 and 4.2]{bdrs}. To see (ix), we have by Lemma \ref{nongg}, (iii) and (iv) that 
$$(-i_X+\frac{r}{2}(n+1-i_X))h_{|F} = (\pi^*(K_X+\det \E))_{|F} = \O_{\P^k}(-k-1)$$
and, since $h_{|F} \cong \O_{\P^k}(1)$, we get that
$$-i_X+\frac{r}{2}(n+1-i_X)=-k-1=-n-r+\nu(\E)$$
that is $\nu(\E)=\frac{(r+2)(n-i_X)+3r}{2}$. Moreover, by the above,
\begin{equation}
\label{ind}
-i_X+\frac{r}{2}(n+1-i_X)=-k-1 \le -2
\end{equation}
hence
\begin{equation}
\label{can}
-K_{\P(\E)}=r\xi+(i_X-\frac{r}{2}(n+1-i_X))h
\end{equation}
is ample by (vii). Also, $-K_{\P(\E)}+\pi^*(K_X)+th=r\xi+(t-\frac{r}{2}(n+1-i_X))h$ is nef non ample by (vii) if and only if $t=\frac{r}{2}(n+1-i_X)$, hence $\tau(\E)=\frac{r}{2}(n+1-i_X)$. Hence (ix) is proved. As for (x), by (ix), the lower bound for $r$ is equivalent to (vi). The upper bound $r \le \frac{2i_X-4}{n+1-i_X}$ follows from \eqref{ind}. As the latter can be rewritten as $\frac{r(n+1)+4}{r+2} \le i_X$, we see that $r \le n-3$ when $X \not\cong Q_n$, since $i_X \le n-1$. When $r=n-3$ we deduce from \eqref{indi} that $i_X \ge n-1$, hence, since $X \not\cong Q_n$, as is well known, $i_X = n-1$ and this proves (x).
\end{proof}

\begin{remark}
\label{spinor}
It follows from Proposition \ref{uno}(vii) that, if $(X, \O_X(1),\E)$ is as in Example \ref{qu}, then $\P(\E)$ is an example of a variety $Y$ with $\Nef(Y)=\Pseff(Y)$.
\end{remark}

\section{Two special cases}

We now prove the main theorem in the case $r=2$ and $b_4(X) = 1$.

\begin{thm}
\label{tre}
Let $X \subseteq \P^N$ be a smooth irreducible variety of dimension $n \ge 2$ such that $\rho(X)=1$ and either $b_4(X) = 1$ or $X$ is Fano with $i_X > \max\{\frac{2n}{3},\frac{n+4}{2}\}$. If $X \not\cong \P^n, X \not\cong Q_n$, any rank $2$ Ulrich bundle on $X$ is big.
\end{thm}
\begin{proof}
Assume that there is a non-big rank $2$ Ulrich bundle $\E$ on $X$. Since $X \not\cong \P^n$ we have that $d \ge 2$. Hence, we know by Proposition \ref{uno}(i)-(ii) that $X$ is Fano, so that $b_2(X)=\rho(X)=1$, and $\Pic(X) \cong \Z \O_X(1)$. If $i_X > \max\{\frac{2n}{3},\frac{n+4}{2}\}$, then $b_4(X) = 1$ by Lemma \ref{b4}. Therefore, we can assume anyway that $b_4(X) = 1$. Note that $\E$ is indecomposable, for otherwise a direct summand of $\E$ must be $\O_X(lH)$, for some $l \in \Z$. Hence, $\O_X(lH)$ is an Ulrich line bundle and this implies that $X \cong \P^n$ (see for example \cite[Lemma 4.2(vi)]{aclr}), a contradiction. 

Now, Proposition \ref{uno}(ix)-(x) gives that $\E$ is Fano with $\tau(\E)=n+1-i_X< i_X-1$. Then \cite[Cor.~6.6]{mos} implies that, if $\E'=\E(-\lfloor \frac{n+2-i_X}{2} \rfloor H)$, then $(X, \E')$ is one of the following:
\begin{itemize}
\item [(1)] $(\P^2, T_{\P^2})$
\item [(2)] $(\P^3,\N)$, $\mathcal N$ a null-correlation bundle
\item [(3)] $(Q_3, \G)$, $\G$ the restriction of a spinor bundle
\item [(4)] $(Q_5, \mathcal C)$, $\mathcal C$ a Cayley bundle
\item [(5)] $(K(G_2), \mathcal Q)$, $\mathcal Q$ the restriction of the universal quotient bundle.
\end{itemize}
Cases (1)-(4) are excluded by hypothesis. In case (5), we have by \cite[Ex.~6.4]{mos} that $n=5$ and $i_X=3$, contradicting Proposition \ref{uno}(x).
\end{proof}

Next, we prove the main theorem for Del Pezzo manifolds.

\begin{thm} 
\label{dp}
Let $X \subset \P^N$ be a Fano manifold of dimension $n \ge 2$ with $\rho(X)=1$ and index $n-1$. Then any Ulrich bundle on $X$ is big.
\end{thm}
\begin{proof}
Assume that there is a non-big rank $r$ Ulrich bundle $\E$ on $X$. By Theorem \ref{nle4} we can assume that $n \ge 5$. By definition, there is an ample line bundle $L$ on $X$ such that $-K_X=(n-1)L$, hence Proposition \ref{uno}(ii) gives that $L=\O_X(1)$. 

The classification of Del Pezzo manifolds (see for example \cite[pages 860-861]{lp}, \cite[Table, page 710]{f}) implies, in the case $\rho(X)=1$, that $(X, \O_X(1))$ is one of the following:
\begin{enumerate}
\item $X \subset \P^{n+1}, n \ge 5$ a smooth cubic hypersurface.
\item $X \subset \P^{n+2}, n \ge 5$ a smooth complete intersection of two quadrics.
\item $X={\mathbb G}(1,4) \subset \P^9$ in the Pl\"ucker embedding, hence $n=6$.
\item $X={\mathbb G}(1,4) \cap H \subset \P^8$ a hyperplane section, hence $n=5$. 
\end{enumerate}

First, suppose that $n=5$. It follows from Proposition \ref{uno}(v) and (x) that $r=2$. In cases (1) and (2), we have that $b_4(X) = 1$ by Lefschetz's hyperplane theorem, and we get a contradiction by Theorem \ref{tre}. In case (4), we have that $c_1(\E)=2H$ by Proposition \ref{uno}(iv) and it follows from \cite[Proof of case (A.3), page 10]{lms1} (or from \cite[Lemma 3.1(viii)]{lr}) that $c_2(\E)H^3=7$. Therefore 
$$s_5(\E^*)=c_1(\E)^5-4c_1(\E)^3c_2(\E)+3c_1(\E)c_2(\E)^2=6c_2(\E)^2H-64 \ne 0$$ 
for otherwise $3$ divides $64$. Therefore $\E$ is big by \cite[Rmk.~2.2]{lm}.

Hence the theorem is established for Del Pezzo manifold of dimension $5$ with $\rho(X)=1$. 

Now assume that $n \ge 6$. Since $r \le n-3$ by Proposition \ref{uno}(x), case (1) is excluded by \cite[Thm.~1]{lrt}, while case (2) is excluded by \cite[Thm.1.1]{es2} (applied to a general hyperplane section when $n$ is even). Thus, only case (3) remains. Hence $X={\mathbb G}(1,4) \subset \P^9$, so that $n=6, i_X=5$ and $r \in \{2, 3\}$ by Proposition \ref{uno}(x). If $r=2$, Proposition \ref{uno}(ix) gives that $\nu(\E)=5$, hence $\E_{|X_5}$ is Ulrich and not big by Lemma \ref{sez}(ii). But $X_5$ is still Del Pezzo, hence this contradicts the case $n=5$ that we have proved above, since $\rho(X_5)=1$ by Lefschetz's theorem. 

Now consider the case $r=3$. We will prove that $\E$ is big by computing 
\begin{equation}
\label{s6}
s_6(\E^*)=c_1(\E)^6-5c_1(\E)^4c_2(\E)+4c_1(\E)^3c_3(\E)+6c_1(\E)^2c_2(\E)^2-6c_1(\E)c_2(\E)c_3(\E)-c_2(\E)^3+c_3(\E)^2.
\end{equation}
We will do calculations in the standard notation of Schubert classes, that we denote by $\sigma_{\lambda}$.
 
We have $H = \sigma_1$ and $c_2(X)=11\sigma_2+12\sigma_{1,1}$ by \cite[Lemma 4.4]{pz}. Then, $c_1(\E)=3H$ by Proposition \ref{uno}(iv), $c_2(\E)\cdot H^4=18$ by \cite[Lemma 3.1(viii)]{lr} and $c_3(\E)\cdot H^3=8$ by \cite[Lemma 5.1(vii)]{blv}. In order to calculate the other quantities, write, for some $a, b, c, d \in \Z$,
$$c_2(\E)=a\sigma_2+b\sigma_{1,1} \ \hbox{and} \ c_3(\E)=c \sigma_3+d\sigma_{2,1}.$$
Using dual Schubert classes, we get 
$$c_2(\E)\sigma_{3,1}=(a\sigma_2+b\sigma_{1,1})\sigma_{3,1}=a \ \hbox{and} \ c_2(\E)\sigma_{2,2}=(a\sigma_2+b\sigma_{1,1})\sigma_{2,2}=b$$
and, since $c_2(\E)$ is represented by an effective $4$-cycle on $X$ (because $\E$ is globally generated), while $\sigma_{3,1}$ and $\sigma_{2,2}$ are nef, we find that
\begin{equation}
\label{a}
a \ge 0, b \ge 0.
\end{equation}
By Pieri's formula we get 
$$H^2=\sigma_2+\sigma_{1,1}, H^3=\sigma_3+2\sigma_{2,1}, H^4=3\sigma_{3,1}+2\sigma_{2,2}.$$
Hence 
$$18=c_2(\E)\cdot H^4=(a\sigma_2+b\sigma_{1,1})(3\sigma_{3,1}+2\sigma_{2,2})=3a+2b$$
which implies, using \eqref{a}, that
\begin{equation}
\label{b}
0 \le b \le 9.
\end{equation}
Also,
$$8=c_3(\E)\cdot H^3=(c \sigma_3+d\sigma_{2,1})(\sigma_3+2\sigma_{2,1})=c+2d.$$ 
Moreover, 
$$c_2(\E)^2\cdot H^2=2a^2+2ab+b^2 \ \hbox{and} \ c_2(\E)^3=a^3+3a^2b+b^3.$$ 
Computing the remaining cases, we obtain
$$c_2(\E)c_3(\E)H=(a\sigma_2+b\sigma_{1,1})(c\sigma_3+d\sigma_{2,1}) \sigma_1=ac+ad+bd$$
and
$$c_3(\E)^2=(c\sigma_3+d\sigma_{2,1})^2=c^2\sigma_3^2+2cd\sigma_3\sigma_{2,1}+d^2\sigma_{2,1}^2=c^2+d^2.$$
From the above, using $H^6 = 5$ and replacing $a=\frac{18-2b}{3}$ and $c=8-2d$ in \eqref{s6}, we obtain
$$s_6(\E^*)=5d^2-30bd+76d-\frac{55}{27}b^3+46b^2-156b+91.$$
Since $0 \le b \le 9$ by \eqref{b}, one easily checks that there exists no integer solution of the equation $s_6(\E^*) = 0$, and therefore $\E$ is big by \cite[Rmk.~2.2]{lm}.

This completes the proof of case (3) and hence of the theorem.
\end{proof}

\section{Proof of Theorem \ref{main}}

In this section we will prove Theorem \ref{main}.

\renewcommand{\proofname}{Proof of Theorem \ref{main}}
\begin{proof}
In the cases (i) and (ii), we have that the Ulrich bundles listed are not big by Examples \ref{sp} and \ref{qu}. 

Vice versa, assume that there is a non-big rank $r$ Ulrich bundle $\E$ on $X$. If either $X \cong \P^n$ or $X \cong Q_n$, then we fall in one of the cases (i)-(ii), by Lemma \ref{pq}.

Therefore, we assume from now on that $X \not\cong \P^n, X \not\cong Q_n$ and we will reach a contradiction. 

We set $c=N-n= \codim_{\P^N}X, d= \deg X$. By Theorem \ref{nle4} and Proposition \ref{uno}(ii), (iii) and (v), we have that 
$$c \ge 1, \ d \ge 3, \ n \ge 5,\ r \ge 2$$
and, for $H \in |\O_X(1)|$, that $\Pic(X) \cong \Z H$ and $-K_X=i_XH$, with $i_X \ge 3$.

We first deal with the case of hypersurfaces.

\renewcommand{\proofname}{Proof}
\begin{claim}
\label{c=1}
The case $c=1$ does not occur.
\end{claim}
\begin{proof}
If $n \in \{5, 6\}$, we have that $r \le 3$ by Proposition \ref{uno}(x), but this is excluded by \cite[Thm.~2]{lr}. If $n \ge 7$, we know that $r \ge \frac{3(n-2)}{4}$ by \cite[Thm.~1]{lrt}. On the other hand, Theorem \ref{dp} implies that $i_X \le n-2$ and then Proposition \ref{uno}(x) gives that
$$r \le \frac{2i_X-4}{n+1-i_X} \le \frac{2(n-4)}{3}$$
contradicting the above lower bound.
\end{proof}
\renewcommand{\proofname}{Proof}
We continue the proof of Theorem \ref{main}. By Claim \ref{c=1} we can assume, from now on, that
$$c \ge 2.$$
We will now divide the proof into several different cases according to the behavior of the map
$$\varphi=\varphi_{\xi} : Y \to Z=\varphi(Y) \subseteq \P H^0(\E)$$ 
where $Y=\P(\E)$ as in Definition-Notation \ref{basica}. We will also use Remark \ref{basic2}. 

We first claim that we can assume that $\dim F=1$ for a general fiber of $\varphi$. 

Indeed, set $\nu=\nu(\E)$, so that $\dim Z = \nu$ and $r-1 \le \nu \le n+r-2$ by Lemma \ref{sez}(ii). If $\nu=n+r-2$, then $\dim Z=n+r-2=\dim Y-1$, hence $\dim F =1$. If $\nu \le n+r-3$, set 
$$N'=N-n+\nu-r+2, \ n'=\nu-r+2$$ 
and consider 
$$X':=X_{\nu-r+2} \subset \P^{N'}.$$
It follows by Lemma \ref{sez}(ii) that 
$$\E':=\E_{|X'}$$ 
is rank $r$ and not big. If either $X' \cong \P^{n'}$ or $X' \cong Q_{n'}$, then, respectively, $\O_{X'}(1) \cong \O_{\P^{n'}}(1)$ or $\O_{X'}(1) \cong \O_{Q_{n'}}(1)$ by Lemma \ref{pq}, but then we get that $X \cong \P^n$ or $X \cong Q_n$, a contradiction. Therefore $X' \not\cong \P^{n'}, X' \not\cong Q_{n'}$ and $n' \ge 5$ by Theorem \ref{nle4}. Now $i_{X'}=i_X-n+\nu-r+2$ and then $n'-i_{X'}=n-i_X$, so that, by Proposition \ref{uno}(ix),
$$\nu(\E')=\frac{(r+2)(n'-i_{X'})+3r}{2}=\frac{(r+2)(n-i_X)+3r}{2}=\nu(\E)=\nu.$$ 
Therefore, if $\varphi': Y' \to Z' \subseteq \P H^0(\E')$, we get, for a general fiber $F'$ that 
$$\dim F'=\dim Y'-\nu(\E')=n'+r-1-\nu=1.$$
Since we are arguing by contradiction, we can replace $X$ with $X'$ and $\E$ with $\E'$. Therefore, as above, we have again $n \ge 5, \Pic(X) \cong \Z H, -K_X=i_XH$, with $i_X \ge 3$ and now the general fiber of the map
$$\varphi=\varphi_{\xi} : Y \to Z=\varphi(Y) \subseteq \P H^0(\E)$$ 
is $\P^1$. Moreover, we have
\begin{equation}
\label{q}
q:=n+1-i_X \ge 3
\end{equation}
by Theorem \ref{dp} and 
$$\dim Z = \nu(\E)=\frac{(r+2)(q-1)+3r}{2}$$
by Proposition \ref{uno}(ix), so that,
$$1=\dim Y - \dim Z=n-q-\frac{rq}{2}.$$
and then 
\begin{equation}
\label{enne}
n=q+1+\frac{rq}{2}, \ i_X=2+\frac{rq}{2}.
\end{equation}
Moreover,
\begin{equation}
\label{c1e}
c_1(\E)=\frac{rq}{2}H
\end{equation} 
by Proposition \ref{uno}(iv) and therefore,
\begin{equation}
\label{ky}
-K_Y=r\xi+2h.
\end{equation} 

We start now to consider several cases according to the dimension of the fibers.

\renewcommand{\proofname}{Proof}
\begin{claim}
\label{tuttep1}
The case $\dim F=1$, for every fiber $F$ of $\varphi$, does not occur.
\end{claim}
\begin{proof}
Assume that $\dim F=1$, for every fiber $F$ of $\varphi$. We have $(F,(\xi+h)_{|F}) \cong (\P^1,\O_{\P^1}(1))$ for all $F$. Since all fibers are connected, in the Stein factorization, we can factor $\varphi : Y \to Z$ through the normalization $\widetilde{\varphi} : Y \to \widetilde Z$ and a bijective morphism $g: \widetilde Z \to Z$. Hence the fibers of $\widetilde{\varphi}$ are the same as the ones of $\varphi$ and, since $\xi+h$ is very ample by Proposition \ref{uno}(viii), we can apply \cite[Prop.~3.2.1]{bs}. It follows that $\widetilde Z$ is smooth and setting $\O_{\widetilde Z}(1)=g^*\O_Z(1)$, we have 
$$\xi=\widetilde{\varphi}^*\O_{\widetilde Z}(1) \ \hbox{and} \ (Y,h) \cong (\P(\F),\O_{\P(\F)}(1))$$
for some rank $2$ bundle $\F$ on $\widetilde Z$. Moreover, the relative tangent bundle $T_{Y/\widetilde Z}$ is a line bundle such that ${T_{Y/\widetilde Z}}_{|F} \cong \O_{\P^1}(2)$ and therefore, since $\Pic(Y)$ is generated by $\xi$ and $h$, there is an integer $m$ such that $T_{Y/\widetilde Z} \cong \O_Y(2h-m\xi)$. Since $T_{Y/\widetilde Z} \cong 2h-\widetilde{\varphi}^*(\det \F)$, it follows that
$\widetilde{\varphi}^*(\det \F) \cong m\xi \cong \widetilde{\varphi}^*\O_{\widetilde Z}(m)$ and therefore, since $\widetilde{\varphi}_*\O_Y \cong \O_{\widetilde Z}$, the projection formula gives 
\begin{equation}
\label{pb}
\det \F \cong \widetilde{\varphi}_*(\widetilde{\varphi}^*(\det \F)) \cong \widetilde{\varphi}_*(\widetilde{\varphi}^*\O_{\widetilde Z}(m)) \cong \O_{\widetilde Z}(m).
\end{equation}
We now show that 
\begin{equation}
\label{m}
m \ge 2.
\end{equation}
Indeed, we have a commutative diagram
$$\xymatrix{0 \ar[r] & \O_Y(2h-m\xi) \ar[r] \ar[dr]^f & T_Y \ar[d]^{d\pi} \\ & & \pi^*T_X}$$ 
with $f \ne 0$, hence injective. Thus, we get an injection $\O_Y \hookrightarrow \pi^*(T_X(-2)) \otimes \O_Y(m\xi)$. Pushing forward to $X$, gives another injection
$$\O_X \cong \pi_*\O_Y \hookrightarrow \pi_*\big(\pi^*(T_X(-2)) \otimes \O_Y(m\xi)\big) \cong T_X(-2) \otimes \pi_*\big(\O_Y(m\xi)\big)$$ 
and therefore
$$H^0(T_X(-2) \otimes \pi_*(\O_Y(m\xi))) \ne 0.$$
Hence, $m \ge 0$ and if $m \in \{0,1\}$, the above contradicts Lemma \ref{ul}(i)-(ii), since $\pi_*\O_Y \cong \O_X$, $\pi_*(\O_Y(\xi)) \cong \E$. Thus $m \ge 2$, proving \eqref{m}. 

Next, observe that Riemann-Roch's theorem gives a polynomial $P(x,y) \in \Q[x,y]$ such that, for any $a,b \in \Z$, we have
$$P(a,b)=\chi(\O_Y(a\xi+bh)).$$
We first claim that
\begin{equation}
\label{rs}
P(\alpha,\beta)=-P(\alpha+m(\beta+1),-\beta-2), \ \hbox{for all} \ \alpha, \beta \in \Q.
\end{equation}
To see the above, it is clearly enough to prove that
\begin{equation}
\label{rs1}
P(a,b)=-P(a+m(b+1),-b-2), \ \hbox{for all} \ a, b \in \Z, b \ge 0.
\end{equation}
To show the latter, fix $a, b \in \Z$ with $b \ge 0$. Observe that
\begin{equation}
\label{rs2}
R^p \widetilde{\varphi}_*\O_Y(a\xi+bh)=0 \ \hbox{for} \ p>0.
\end{equation}
In fact, since $\O_Y(a\xi+bh) \cong \O_Y(bh) \otimes \widetilde{\varphi}^*\O_{\widetilde Z}(a)$, \eqref{rs2} follows because $R^p \widetilde{\varphi}_*\O_Y(bh)=0$ for $p>0$ by \cite[Ex.~III.8.4(a)]{h1}. Then, $H^i(\O_Y(a\xi+bh)) \cong H^i(\widetilde{\varphi}_*\O_Y(a\xi+bh))$ and the Leray's spectral sequence gives that
$$P(a,b)=\chi(\O_Y(a\xi+bh))=\chi(\widetilde{\varphi}_*\O_Y(a\xi+bh)).$$
On the other hand, since $\F$ has rank $2$, we have by \eqref{pb} that $\F^* \cong \F(-\det \F) \cong \F \otimes \O_{\widetilde Z}(-m)$, hence $(S^b \F^*)(-\det \F) \cong (S^b \F) \otimes \O_{\widetilde Z}(-m(b+1))$. Now, \cite[Ex.~III.8.4(c)]{h1} gives that 
$$R^1 \widetilde{\varphi}_*\O_Y(-(b+2)h) \cong (S^b \F^*)(-\det \F) \cong (S^b\F) \otimes \O_{\widetilde Z}(-m(b+1))$$
and therefore
\begin{equation}
\label{rs3}
\begin{aligned}
R^1 \widetilde{\varphi}_*\O_Y((a+m(b+1))\xi-(b+2)h) & \cong R^1 \widetilde{\varphi}_*\O_Y(-(b+2)h) \otimes \O_{\widetilde Z}(a+m(b+1)) \\ & \cong (S^b \F) \otimes \O_{\widetilde Z}(a) \cong \widetilde{\varphi}_*\O_Y(a\xi+bh)
\end{aligned}
\end{equation}
Moreover, since $-b-2 < 0$, we have that $R^p \widetilde{\varphi}_*\O_Y(-(b+2)h)=0$ for $p \ne 1$ and therefore also 
$$R^p \widetilde{\varphi}_*\O_Y((a+m(b+1))\xi-(b+2)h)=0 \ \hbox{for} \ p \ne 1.$$
It follows from the invariance of the Euler characteristic in Leray's spectral sequence \cite[Lemma 33.33.5, \href{https://stacks.math.columbia.edu/tag/0BEK}{Tag 0BEK}]{st} and \eqref{rs3} that
$$\begin{aligned}
P(a+&m(b+1),-b-2) = \chi(\O_Y((a+m(b+1))\xi-(b+2)h))= \\ 
& \sum\limits_{p \ge 0}(-1)^p\chi(R^p \widetilde{\varphi}_*\O_Y((a+m(b+1))\xi-(b+2)h))=-\chi(R^1 \widetilde{\varphi}_*\O_Y((a+m(b+1))\xi-(b+2)h)) = \\ 
& -\chi(\widetilde{\varphi}_*\O_Y(a\xi+bh)) = -\chi(\O_Y(a\xi+bh))=-P(a,b)
\end{aligned}$$
and \eqref{rs1} is proved. Therefore, so is \eqref{rs}.

Observe now that $P(-j,b)=0$ for $j \in \{1,\ldots,r-1\}$ and all $b \in \Z$: Indeed, $R^p \pi_*(-j\xi)=0$ for $p \ge 0$ by \cite[Ex.~III.8.4(a)]{h1}, hence $H^i(\O_Y(-j\xi+bh))=0$ for $i \ge 0$ and therefore $P(-j,b)=\chi(\O_Y(-j\xi+bh))=0$. It follows that
\begin{equation}
\label{rs4}
P(-j,\beta)=0 \ \hbox{for} \ j \in \{1,\ldots,r-1\}, \beta \in \Q.
\end{equation}
Moreover, if $b \in \Z$, since $R^p \pi_*(\xi)=0$ for $p>0$ by \cite[Ex.~III.8.4(a)]{h1}, we have
$$P(1,b)=\chi(\O_Y(\xi+bh))=\chi(\pi_*(\xi)(b))=\chi(\E(b))=rd\binom{b+n}{n}=\frac{1}{n!}rd(b+n)(b+n-1)\cdots (b+1)$$
by \cite[Lemma 3.1(ix)]{lr}. Since the above holds for every $b \in \Z$, we deduce that
\begin{equation}
\label{rs5}
P(1,\beta)=\frac{1}{n!}rd(\beta+n)(\beta+n-1)\cdots (\beta+1)
\end{equation}
for any $\beta \in \Q$. Also observe that
\begin{equation}
\label{rs6}
P(1,\beta)=0 \iff \beta \in \{-1,\ldots,-n\}.
\end{equation}
Assume that $r \ge 3$. For $j \in \{1,2\}$, let 
$$\beta=-1+\frac{j+1}{m}$$ 
so that $-j+m(\beta+1)=1$. Then, using \eqref{rs4} and \eqref{rs}, we have that
$$0=P(-j,\beta)=-P(1,-\beta-2)$$
and therefore, by \eqref{rs6}, there is an $i \in \{1,\ldots,n\}$ with $\beta=i-2 \in \Z$, that is $m \mid j+1$. But then $m \mid 2$ and $m \mid 3$, a contradiction since $m \ge 2$ by \eqref{m}.

Hence $r=2$. Then, \cite[Thm.~1.1]{wa} implies that $X \cong K(G_2)$. Since $b_4(K(G_2))=1$ (see for example \cite[Ex.~6.4]{mos}), this case is excluded by Theorem \ref{tre}. This completes the proof of Claim \ref{tuttep1}.
\end{proof}
\renewcommand{\proofname}{Proof}
We continue the proof of Theorem \ref{main}. 

It follows by Claim \ref{tuttep1} that there is at least one fiber of $\varphi$ with dimension at least $2$.

Let $F$ be a fiber with $\dim F \ge 2$. Restricting to $F$ the exact sequence
$$0 \to \Omega_{Y/X} \to (\pi^* \E)(-\xi) \to \O_Y \to 0$$
we get, with the usual identification $F \cong \pi(F) \subset X$,
$$0 \to (\Omega_{Y/X})_{|F} \to \E_{|F} \to \O_F \to 0.$$ 
Since $\E_{|F}$ is globally generated, because $\E$ is, the above sequence splits and, considering the rank $r-1$ bundle $\G=(\Omega_{Y/X})_{|F}$, we find that
\begin{equation}
\label{g}
\E_{|F} \cong \O_F \oplus \G
\end{equation}
and \eqref{g}, \eqref{c1e} give that 
\begin{equation}
\label{c1g}
c_1(\G)=\frac{rq}{2}.
\end{equation} 
Now the diagram
$$\xymatrix{& & 0 \ar[d] & 0 \ar[d] & \\ 0 \ar[r] & N_{\pi(F)/X}^* \ar[r] & (\pi^*\Omega_X)_{|F} \cong {\Omega_X}_{|\pi(F)} \ar[d] \ar[r] & \Omega_{\pi(F)} \ar[r] \ar[d]^{\cong} & 0 & \\ 0 \ar[r] & N_{F/Y}^* \ar[r] & {\Omega_Y}_{|F} \ar[r] \ar[d] & \Omega_F \ar[r]  \ar[d] & 0 & \\ & & \G \ar[d] & 0 & \\ & & 0 & &}$$
gives, by snake lemma, with the usual identification $F \cong \pi(F) \subset X$, the exact sequence
\begin{equation}
\label{g1}
0 \to N_{F/X}^* \to N_{F/Y}^* \to \G \to 0.
\end{equation}
Consider the Ulrich dual, which is also Ulrich, $\E^*(K_X+(n+1)H)=\E^*(q)$ and the rank $r-1$ bundle 
\begin{equation}
\label{v1}
\V=\G^*(q).
\end{equation}
We have, by \eqref{g},
\begin{equation}
\label{c1v}
c_1(\V)=c_1(\G^*(q))=-c_1(\G(-q))=-c_1(\G)+(r-1)q=\frac{(r-2)q}{2}.
\end{equation}
Let $P \subseteq F$ be a linear plane and set $\V'=\V_{|P}$. We have, by \eqref{g},
\begin{equation}
\label{nu1}
\E^*(q)_{|P} \cong \O_P(q) \oplus \V'.
\end{equation}
We claim that
\begin{equation}
\label{nu2}
h^0(\V'(-1))=\frac{1}{8}q(r-2)(rq-2q+2)-c_2(\V').
\end{equation}
To see the above, since $\E^*(q)$ is Ulrich, we can apply Lemma \ref{ul}(iii) to it, and get $H^1(\E^*(q)_{|P}(-1))=H^2(\E^*(q)_{|P}(-1))=0$ and therefore also $H^1(\V'(-1))=H^2(\V'(-1))=0$ by \eqref{nu1}. Hence, to prove \eqref{nu2}, it remains to show that 
\begin{equation}
\label{nu3}
\chi(\V'(-1))=\frac{1}{8}q(r-2)(rq-2q+2)-c_2(\V').
\end{equation}
In order to prove the above, using \eqref{c1v}, we find
$$c_1(\V'(-1))=\frac{(r-2)q}{2}-r+1$$
and
$$c_2(\V'(-1))=c_2(\V')-\frac{(r-2)^2q}{2}+\binom{r-1}{2}=c_2(\V')-\frac{1}{2}(r-2)(rq-2q-r+1)$$
and applying Riemann-Roch we get \eqref{nu3}, hence also \eqref{nu2}. 

We now deal with the next case.

\renewcommand{\proofname}{Proof}
\begin{claim}
\label{nop2}
There is no fiber $F$ of $\varphi$ with $\dim F=2$.
\end{claim}
\begin{proof}
Assume that there is a fiber $F$ of $\varphi$ with $\dim F=2$ and let $P=F$ and $\V'=\V$. We will prove that
\begin{equation}
\label{s1}
\sigma_2(N_{P/X}^*(1))=-h^0(\V(-1)).
\end{equation}
The above implies, by Lemma \ref{plane} and \eqref{q}, the contradiction
$$2 \le q-1=n-i_X \le \sigma_2(N_{P/X}^*(1)) \le 0$$
and therefore Claim \ref{nop2} will be proved once shown \eqref{s1}. To see the latter, we have by \eqref{ky} that
$$c_1(N_{P/Y}^*)=c_1(P)-c_1(Y)_{|P}=1$$
and therefore
$$c_1(N_{P/Y}^*(1))=n+r-2.$$
Since $N_{P/Y}^*$ is globally generated and non-big by Lemma \ref{nf}, we have by \cite[Rmk.~2.2]{lm} that 
$$0=\sigma_2(N_{P/Y}^*)=1-c_2(N_{P/Y}^*).$$ 
Hence $c_2(N_{P/Y}^*)=1$ and then
\begin{equation}
\label{c2n}
c_2(N_{P/Y}^*(1))=c_2(N_{P/Y}^*)+(n+r-4)c_1(N_{P/Y}^*)+\binom{n+r-3}{2}=n+r-3+\binom{n+r-3}{2}=\binom{n+r-2}{2}.
\end{equation}
To compute $c_2(N_{P/X}^*(1))$, we will use \eqref{g1} twisted by $\O_P(1)$, that is
\begin{equation}
\label{ns}
0 \to N_{P/X}^*(1) \to N_{P/Y}^*(1) \to \G(1) \to 0.
\end{equation}
Since, by \eqref{c1g}, 
\begin{equation}
\label{c1g1}
c_1(\G(1))=c_1(\G)+r-1=\frac{rq}{2}+r-1
\end{equation} 
we have, by \eqref{v1} and \eqref{c1v}, that
\begin{equation}
\label{c2g1}
c_2(\G(1))=c_2(\V^*(q+1))=c_2(\V(-q-1))=c_2(\V)-\frac{(r-2)^2q(q+1)}{2}+(q+1)^2\binom{r-1}{2}.
\end{equation}
Then, using \eqref{c2n}, \eqref{ns} and Lemma \ref{plane}, we get
$$c_2(N_{P/X}^*(1))=\binom{n+r-2}{2}-c_2(\G(1))-q(\frac{rq}{2}+r-1)$$
and combining with \eqref{c2g1} we find, using also \eqref{enne} and \eqref{nu2},
$$\begin{aligned}
\sigma_2(N_{P/X}^*(1)) & =q^2-\binom{n+r-2}{2}+c_2(\V)-\frac{(r-2)^2q(q+1)}{2}+(q+1)^2\binom{r-1}{2}+q(\frac{rq}{2}+r-1) \\ & = -\frac{1}{8}q(r-2)(rq-2q+2)+c_2(\V)=-h^0(\V(-1))
\end{aligned}$$
so that \eqref{s1} follows. This completes the proof of Claim \ref{nop2}.
\end{proof}
\renewcommand{\proofname}{Proof}

We continue the proof of Theorem \ref{main}. We exclude $r \ge 4$.

\renewcommand{\proofname}{Proof}
\begin{claim}
\label{r<4}
The case $r \ge 4$ does not occur.
\end{claim}
\begin{proof}
Since $r \ge 4$, we have by \eqref{enne} and \eqref{q} that $i_X=2+\frac{rq}{2}>\max\{\frac{2n}{3},\frac{n+4}{2}\}$ and then $H^4(X;\Z) \cong \Z H^2$ by Lemma \ref{b4}. Therefore there are $a, b \in \Z$ such that
\begin{equation}
\label{c2ab}
c_2(X)=a H^2, c_2(\E)=b H^2.
\end{equation}
By Claims \ref{tuttep1} and \ref{nop2} there is a fiber $F$ of $\varphi$ with $\dim F \ge 3$. Choose a linear plane $P \subset F$ and set
\begin{equation}
\label{ep}
\sigma=\sigma_2(N_{P/X}^*(1)).
\end{equation}
We will first prove that
\begin{equation}
\label{gs}
a=\frac{1}{8}(q^2r^2+4q^2+6qr+4q+8)-\sigma.
\end{equation}
We have
$$c_2(T_X(-1))=c_2(X)-(n-1)c_1(X)H+\binom{n}{2}H^2=(a-(n-1)i_X+\binom{n}{2})H^2$$
and therefore
\begin{equation}
\label{c2a}
c_2(T_X(-1)_{|P})=a-(n-1)i_X+\binom{n}{2}.
\end{equation}
Since $c_1(T_P(-1))=c_2(T_P(-1))=1$ by the Euler sequence and $c_1(N_{P/X}(-1))=-q$ by Lemma \ref{plane}, we get, from the exact sequence 
$$0 \to T_P(-1) \to {T_X(-1)}_{|P} \to N_{P/X}(-1) \to 0$$
using \eqref{ep} and \eqref{c2a}, that
$$\begin{aligned}
a-(n-1)i_X+\binom{n}{2} & =c_2(T_X(-1)_{|P}) \\
& =1+c_2(N_{P/X}(-1))-q \\
& = 1-q+c_1(N_{P/X}(-1))^2-\sigma_2(N_{P/X}^*(1)) \\
& = 1-q+q^2-\sigma
\end{aligned}$$
hence
$$a=(n-1)i_X-\binom{n}{2}+1-q+q^2-\sigma$$
and, using \eqref{enne}, this proves \eqref{gs}. Next, we claim that 
\begin{equation}
\label{c2e}
b=\frac{r}{24}(3q^2r-2q^2+2q-2\sigma).
\end{equation}
In fact, since $d=H^n$, \cite[Lemma 3.1(viii)]{lr} gives, using \eqref{c1e},
$$\begin{aligned}
b d & = c_2(\E) H^{n-2}=\frac{1}{2}\big(c_1(\E)^2-c_1(\E) K_X\big) H^{n-2}+\frac{r}{12}\big(K_X^2+c_2(X)-\frac{3n^2+5n+2}{2}H^2\big) H^{n-2} \\
& = d \left(\frac{r^2q^2}{8}+\frac{rqi_X}{4}+\frac{r}{12}\big(i_X^2+a-\frac{3n^2+5n+2}{2}\big)\right)
\end{aligned}$$
that is, using \eqref{enne} and \eqref{gs}, exactly \eqref{c2e}. Now, we have, using \eqref{c1e},
$$c_2(\E^*(q)_{|P})=c_2(\E^*_{|P})+(r-1)c_1(\E^*_{|P})q+\binom{r}{2}q^2=b$$
From \eqref{nu1}, we have $\E^*(q)_{|P} \cong \O_P(q) \oplus \V'$, hence, since $c_1(\V')=\frac{(r-2)q}{2}$ by \eqref{c1v}, we get that
$$b=c_2(\E^*(q)_{|P})=c_2(\V')+\frac{(r-2)q^2}{2}$$
and then
$$c_2(\V')=b-\frac{(r-2)q^2}{2}.$$
Replacing in the above the expression of $b$ in \eqref{c2e}, we find that
\begin{equation}
\label{c2v}
c_2(\V')=\frac{1}{24}(24q^2-2r\sigma+2qr-14rq^2+3q^2r^2).
\end{equation}
Next, set 
$$z=h^0(\V'(-1))$$
so that
$$z=\frac{1}{8}q(r-2)(rq-2q+2)-c_2(\V')$$
by \eqref{nu2}. Replacing in the above the expression of $c_2(\V')$ in \eqref{c2v}, we find the identity
\begin{equation}
\label{ze}
12z=r\big(\sigma+q(q+2)\big)-6q(q+1).
\end{equation}
We will prove that the above identity leads to a contradiction. To this end, we claim that there is an integer $m$ such that
\begin{equation}
\label{zm}
2m(z+\sigma+1)=rq, \ m \ge 2.
\end{equation}
Granted the above for a moment, let us get a contradiction. Recall that $z \ge 0, q \ge 3, r \ge 4$ and $\sigma \ge q-1$ by Lemma \ref{plane}. Then \eqref{zm}  gives
$$3rq=6m(z+\sigma+1)\ge 12(z+\sigma+1)=r(\sigma+q(q+2))-6q(q+1)+12\sigma+12$$
so that
$$(12+r)(q-1) \le \sigma(12+r) \le 3rq-rq(q+2)+6q(q+1)-12$$
and then 
$$(q-1)\big(r(q+1)-6q) \le 0.$$
It follows that $r \le \frac{6q}{q+1}<6$. Moreover, from \eqref{zm} we have that $rq=2m(z+\sigma+1)\ge 2mq$, so that $m \le \frac{r}{2}<3$. Therefore $m=2$ and $r \in \{4,5\}$. If $r=4$, solving for $z$ the system given by \eqref{ze} and \eqref{zm}, we get $8z=-(q-1)(q-2)<0$, a contradiction. If $r=5$, we have $5 \le \frac{6q}{q+1}$, so that $q \ge 5$. Moreover, $4 \mid q$ by \eqref{zm}, hence $q \ge 8$. Solving for $z$ the system given by \eqref{ze} and \eqref{zm}, we get
$$17z=-q^2+\frac{41q}{4}-5<0$$
a contradiction if $q \ge 12$. Hence $q=8$ and the above is $17z=13$, again a contradiction. To finish the proof of Claim \ref{r<4} it remains to prove \eqref{zm}. In order to do that, we return to the map
$$\varphi=\varphi_{\xi} : Y \to Z=\varphi(Y) \subseteq \P H^0(\E)$$ 
where $Y=\P(\E)$ as in Definition-Notation \ref{basica}. We will also use Remark \ref{basic2}. 
 
Note that $\xi$ and $h$ are nef. Recall that the general fiber $L$ of $\varphi$ is $\P^1$ and
\begin{equation}
\label{y1}
\dim Y = n+r-1, \ \dim Z = n+r-2, \ \xi \cdot L=0, \ h \cdot L=1, \ \xi^{n+r-1}=0, \ -K_Y=r\xi+2h.
\end{equation} 
Then
\begin{equation}
\label{int01}
h\xi^{n+r-2}=h \varphi^*\O_Z(1)^{n+r-2}=(\deg Z)L \cdot h= \deg Z >0.
\end{equation}
We now compute some intersections on $Y$. We first claim that
\begin{equation}
\label{int11}
\ch_2(Y) \xi^{n+r-3}=0.
\end{equation}
To see the above, observe first that
$$H^i(t\xi-h)=0 \ \hbox{for} \ t \ge 0, i \ge 0.$$
Indeed, if $i>0$, the vanishing holds by Kodaira's theorem, since $t\xi-h=K_Y+(t+r)\xi+h$ and $(t+r)\xi+h$ is ample by Proposition \ref{uno}(vii). If $i=0$, since $(t\xi-h) \cdot L = -1$, we get that $H^0(t\xi-h)=0$, because a global section vanishes on a general fiber. Thus $P(t)=\chi(t\xi-h)$ is the zero polynomial. By Riemann-Roch, its coefficient of $t^{n+r-3}$, multiplied by $(n+r-3)!$, is
$$0=\left(\frac{h^2}{2}+\frac{hK_Y}{2}+\frac{K_Y^2+c_2(Y)}{12}\right)\xi^{n+r-3}$$
so that \eqref{y1} gives
$$0=\bigl(c_2(Y)-2r\xi h-2h^2\bigr)\xi^{n+r-3}$$
hence
$$c_2(Y)\xi^{n+r-3}=\frac12K_Y^2\xi^{n+r-3}$$
and then
$$\ch_2(Y)\xi^{n+r-3}=\frac{1}{2}\big(K_Y^2-2c_2(Y)\big)\xi^{n+r-3}=0$$
that is \eqref{int11}. Next, we show that there exists an integer $m$ such that
\begin{equation}
\label{int21}
m \ge 2 \ \ \hbox{and} \ \ h^2\xi^{n+r-3}=m h\xi^{n+r-2}.
\end{equation}
Replacing $\varphi$ by the connected-fiber morphism in its Stein factorization, we retain the notation $\varphi:Y\to Z$ and denote by $\O_Z(1)$ the pullback of the original hyperplane bundle. Then $Z$ is normal and the fibers are unchanged scheme-theoretically, since the original fibers are connected and reduced. Let $U$ be the non-empty open subset of $Z \setminus \Sing(Z)$ where the fibers of $\varphi$ are $1$-dimensional. We have that
$$\codim_Z (Z \setminus U) \ge 2$$
for otherwise $\dim (Z \setminus U) \ge \dim Z-1$, but then, being the fibers of $\varphi$ over general points of the components of $Z \setminus U$ of dimension at least $2$, we would have that
$$\dim \varphi^{-1}(Z \setminus U) \ge \dim (Z \setminus U)+2 \ge \dim Z+1=\dim Y$$
hence $\varphi^{-1}(Z \setminus U)=Y$, a contradiction. Now, $\varphi^{-1}(U) \to U$ is a smooth $\P^1$-bundle and, since $h$ has degree one on every fiber, it is a projective bundle. Extend $T_{\varphi^{-1}(U)/U}$ to a coherent subsheaf of $T_Y$, keeping the same notation. We have an inclusion $T_{\varphi^{-1}(U)/U} \subset T_Y$ that defines an exact sequence
$$0 \to T_{\varphi^{-1}(U)/U} \to T_Y \to \G \to 0$$
where $\G$ is the quotient. Observe that, $\G$ is torsion free at the points $y \in \varphi^{-1}(U)$, since
$$\G_y \cong T_{Y,y}/(T_{\varphi^{-1}(U)/U})_y \cong \Im d \varphi_y.$$
Consider the commutative diagram
$$\xymatrix{& 0 \ar[d] & 0 \ar[d] & 0 \ar[d] & \\ 0 \ar[r] & T_{\varphi^{-1}(U)/U} \ar[d]^{\Id} \ar[r] & \Ker(g \circ f) \ar[d] \ar[r] & {\mathcal T}(\G) \ar[r] \ar[d] & 0 & \\ 0 \ar[r] & T_{\varphi^{-1}(U)/U} \ar[r] \ar[d] & T_Y \ar[r]^{f} \ar[d]^{g \circ f} & \G \ar[r] \ar[d]^{g} & 0 & \\ & 0 \ar[r] & \G/{\mathcal T}(\G) \ar[d] \ar[r]^{\Id} & \G/{\mathcal T}(\G) \ar[r] \ar[d] & 0. \\ & & 0 & 0 & }$$
The subsheaf $\F:=\Ker(g \circ f)$ is rank $1$ and reflexive by \cite[Prop.~1.1]{h3}, since it sits in the exact sequence given by the second column of the diagram. Hence $\F$ is a line bundle on $Y$ by \cite[Prop.~1.9]{h3}. Therefore, we can write $\F=\O_Y(sh-m\xi)$ for some $s, m \in\Z$. Moreover, since $\G$ is torsion free at the points $y \in \varphi^{-1}(U)$, on a general fiber $L \subset \varphi^{-1}(U)$, we have that $\F_{|L} \cong (T_{\varphi^{-1}(U)/U})_{|L} \cong \O_{\P^1}(2)$, and then $\F=\O_Y(2h-m\xi)$. Take a general complete intersection curve $C = \cap H_1 \cap \dots \cap H_{n+r-3}$, where the $H_j$'s are general members of $|\O_Z(1)|$. Since $\codim_Z (Z \setminus U) \ge 2$, we have that $C \subset U$ and then the surface $S=\varphi^{-1}(C)$ is a ruled surface over $C$. Moreover, by construction, $\F_{|S}=T_{S/C}$. Since $c_1(T_{S/C})^2=0$, we have that $(2h-m\xi)_{|S}^2=0$ so that
$$(2h-m\xi)^2\xi^{n+r-3}=0.$$
Together with $\xi^{n+r-1}=0$ by \eqref{y1}, this gives 
\begin{equation}
\label{int31}
h^2\xi^{n+r-3}=m h\xi^{n+r-2}
\end{equation} 
proving the second condition in \eqref{int21}. To show that $m \ge 2$, note that the differential of $\pi$ gives a non-zero, hence injective, map $\O_Y(2h-m\xi) \hookrightarrow \pi^*T_X$, hence also an injection $\O_Y \hookrightarrow \pi^*(T_X(-2)) \otimes \O_Y(m\xi)$. Pushing forward to $X$, we get another injection
\begin{equation}
\label{int41}
\O_X \cong \pi_*\O_Y \hookrightarrow \pi_*\big(\pi^*(T_X(-2)) \otimes \O_Y(m\xi)\big) \cong T_X(-2) \otimes \pi_*\O_Y(m\xi).
\end{equation} 
Hence, $m \ge 0$ and if $m \in \{0,1\}$, the above contradicts Lemma \ref{ul}(i)-(ii), since $\pi_*\O_Y \cong \O_X$, $\pi_*(\O_Y(\xi)) \cong \E$. Thus $m \ge 2$, proving \eqref{int21}. We will use the latter to prove \eqref{zm}. To this end, we  compute $\ch_2(Y)=\frac{1}{2}c_1(Y)^2-c_2(Y)$. We know that 
$$c_1(Y)=-K_Y=r\xi+2h$$
by \eqref{y1}. From the sequence of tangent bundles
$$0 \to T_{Y/X} \to T_Y \to \pi^*T_X \to 0$$
and the relative Euler sequence
$$0 \to \O_Y \to (\pi^*\E^*)(\xi) \to T_{Y/X} \to 0$$
we get, using \eqref{c2ab} and \eqref{c1e}, that
$$c_2(Y)=\binom{r}{2}\xi^2+\frac{r(q+4)}{2}\xi h+\big(a+b-\frac{rq(qr+4)}{4}\big)h^2$$
and this gives
$$\ch_2(Y)=\frac{r}{2} \xi^2-\frac{rq}{2}\xi h+(\frac{r^2q^2}{4}+rq+2-a-b)h^2.$$
Now, replacing into the above the expressions of $a$ and $b$ in \eqref{gs} and \eqref{c2e}, we get that
$$\ch_2(Y)=\frac{r}{2} \xi^2-\frac{rq}{2}\xi h+\frac{1}{12}(12+12\sigma-6q-6q^2+r\sigma+2qr+rq^2)h^2$$
and using \eqref{ze}, we find
$$\ch_2(Y)=\frac{r}{2} \xi^2-\frac{rq}{2}\xi h+(z+\sigma+1)h^2.$$
Intersecting with $\xi^{n+r-3}$ and using \eqref{int11}, \eqref{y1} and \eqref{int21} we get
$$0=\ch_2(Y)\xi^{n+r-3}=-\frac{rq}{2}h\xi^{n+r-2}+(z+\sigma+1)h^2\xi^{n+r-3}=\big(-\frac{rq}{2}+(z+\sigma+1)m\big)h\xi^{n+r-2}.$$
Since $h\xi^{n+r-2}>0$ by \eqref{int01}, we deduce that
$$-\frac{rq}{2}+(z+\sigma+1)m=0$$
that is \eqref{zm}. This completes the proof of Claim \ref{r<4}.
\end{proof}
\renewcommand{\proofname}{Proof}

We now complete the proof of Theorem \ref{main}. 

By Claims \ref{tuttep1} and \ref{nop2} there is a fiber $F \cong \P^s$ of $\varphi$ with $s \ge 3$. Moreover, we have by Claim \ref{r<4} that $2 \le r \le 3$. Consider the exact sequence \eqref{g1}
$$0 \to N_{F/X}^* \to N_{F/Y}^* \to \G \to 0$$
and set $\F'=N_{F/X}^*(1), \F=N_{F/Y}^*$ and $\F''=\G$. It follows from Lemma \ref{nf} that $\F$ is globally generated and non-big. Also, by \eqref{ky},
$$c_1(\F)=c_1(N_{F/Y}^*)=-K_F+{K_Y}_{|F}=s-1$$
and
$$c_1(N_{F/X}^*)=-K_F+{K_X}_{|F}=s+1-i_X$$
so that $c_1(\F')=c_1(N_{F/X}^*(1))=1-i_X+n=q \ge 2$ by \eqref{q}. If $r=2$ we have, by \eqref{c1g}, that $\G$ is a line bundle with $c_1(\G)=q$, hence $\F''= \G \cong \O_{\P^s}(q)$. If $r=3$, then $\G$ has rank $2$ and \eqref{c1g} implies that $q$ is even, say $q=2t$, so that $c_1(\G)=3t$. If we set $V=\V=\G^*(2t)$ as in \eqref{v1}, we have that $V$ has rank $2$, $c_1(V)=t$ by \eqref{c1g} and $V(t) \cong \G^*(3t) \cong \G=\F''$. Moreover, as in \eqref{nu1}, we have that $\E^*(q)_{|F} \cong \O_F(q) \oplus \V$, hence $V=\V$ is globally generated since $\E^*(q)$ is Ulrich (being the Ulrich dual of $\E$). Therefore the hypotheses of Lemma \ref{bnogg} are satisfied and the lemma implies that $\F'=N_{F/X}^*(1)$ is not nef. But, as is well known, the conormal bundle sequence of $F \subset X \subset \P^N$ gives a surjection
$$\O_F^{\oplus (N-s)} \twoheadrightarrow N_{F/X}^*(1)$$
hence $N_{F/X}^*(1)$ is globally generated, a contradiction.

This completes the proof of Theorem \ref{main}.
\end{proof}
\renewcommand{\proofname}{Proof}


\begin{thebibliography}{ACLR}

\bibitem[ACLR]{aclr} V.~Antonelli, G.~Casnati, A.~F.~Lopez, D.~Raychaudhury.
\textit{On varieties with Ulrich twisted conormal bundles}.
Proc. Amer. Math. Soc. \textbf{152} (2024), no.~11, 4645-4658.

\bibitem[Ba]{ba} C.~B\u anic\u a.
\textit{Smooth reflexive sheaves}. 
Proceedings of the Colloquium on Complex Analysis and the Sixth Romanian-Finnish Seminar,
Rev. Roumaine Math. Pures Appl. \textbf{36} (1991), no.~9-10, 571-593.

\bibitem[Be]{b} A.~Beauville.
\textit{An introduction to Ulrich bundles}. 
Eur. J. Math. \textbf{4} (2018), no.~1, 26-36.

\bibitem[BDRS]{bdrs} M.~C.~Beltrametti, S.~Di Rocco, A.~J.~Sommese.
\textit{On generation of jets for vector bundles}. 
Rev. Mat. Complut. \textbf{12} (1999), no.~1, 27-45.

\bibitem[BLV]{blv} V.~Buttinelli, A.~F.~Lopez, R.~Vacca.
\textit{Ulrich bundles with $c_2(\E)^2=0$ and connectedness of Ulrich subvarieties}. 
{\it Preprint 2026 (arXiv:2604.22978)}.

\bibitem[BS]{bs} M.~C.~Beltrametti, A.~J.~Sommese.
\textit{The adjunction theory of complex projective varieties}. 
De Gruyter Expositions in Mathematics, \textbf{16}. Walter de Gruyter \& Co., Berlin, 1995.

\bibitem[CH]{ch} M.~Casanellas, R.~Hartshorne.
\textit{Stable Ulrich bundles}. With an appendix by F.~Geiss, F.-O.~Schreyer.
Internat. J. Math. \textbf{23} (2012), no.~8, 1250083, 50 pp.

\bibitem[DPS]{dps} J.~P.~Demailly, T.~Peternell, M.~Schneider.
\textit{Compact complex manifolds with numerically effective tangent bundles}. 
J. Algebraic Geom. \textbf{3} (1994), no. 2, 295-345. 

\bibitem[E]{e} P.~Ellia.
\textit{Chern classes of rank two globally generated vector bundles on $\P^2$}.
Atti Accad. Naz. Lincei Rend. Lincei Mat. Appl. \textbf{24} (2013), no.~2, 147-163.

\bibitem[EH]{eh} D.~Eisenbud, J.~Harris.
\textit{3264 and all that - a second course in algebraic geometry}.
Cambridge University Press, Cambridge, 2016. xiv+616 pp. 

\bibitem[ES1]{es} D.~Eisenbud, F.-O.~Schreyer.
\textit{Resultants and Chow forms via exterior syzygies}. 
J. Amer. Math. Soc. \textbf{16} (2003), no.~3, 537-579.

\bibitem[ES2]{es2} D.~Eisenbud, F.-O.~Schreyer.
\textit{Hyperelliptic curves and Ulrich sheaves on the complete intersection of two quadrics}. 
\'Epijournal G\'eom. Alg\'ebrique [2023-2025], Special volume in honour of Claire Voisin, Art. 21, 25 pp.

\bibitem[F]{f} T.~Fujita.
\textit{On the structure of polarized manifolds with total deficiency one. I}. 
J. Math. Soc. Japan \textbf{32} (1980), no.~4, 709-725.

\bibitem[H1]{h1} R.~Hartshorne.
\textit{Algebraic geometry}. 
Graduate Texts in Mathematics, No. 52. Springer-Verlag, New York-Heidelberg, 1977.

\bibitem[H2]{h2} R.~Hartshorne.
\textit{Varieties of small codimension in projective space}. 
Bull. Amer. Math. Soc. \textbf{80} (1974), 1017-1032. 

\bibitem[H3]{h3} R.~Hartshorne.
\textit{Stable reflexive sheaves}. 
Math. Ann. \textbf{254} (1980), no.~2, 121-176.

\bibitem[IR]{ir2} P.~Ionescu, F.~Russo.
\textit{Remarks on defective Fano manifolds}.
Ann. Univ. Ferrara Sez. VII Sci. Mat. \textbf{63} (2017), no.~1, 133-146.

\bibitem[La1]{la} R.~Lazarsfeld. \textit{Positivity in algebraic geometry, I. Classical Setting: Line Bundles and Linear Series}.
Ergebnisse der Mathematik und ihrer Grenzgebiete. 3. Folge. A Series of Modern Surveys in Mathematics, \textbf{48}. Springer-Verlag, Berlin, 2004.

\bibitem[La2]{laz2} R.~Lazarsfeld. \textit{Positivity in algebraic geometry. II. Positivity for vector bundles, and multiplier ideals}.
Ergebnisse der Mathematik und ihrer Grenzgebiete. 3. Folge. A Series of Modern Surveys in Mathematics, \textbf{49}. Springer-Verlag, Berlin, 2004.

\bibitem[Lo]{lo} A.~F.~Lopez. 
\textit{On the positivity of the first Chern class of an Ulrich vector bundle}.
Commun. Contemp. Math. \textbf{24} (2022), no.~9, Paper No.~2150071, 22 pp.

\bibitem[LM]{lm} A.~F.~Lopez, R.~Mu\~{n}oz. 
\textit{On the classification of non-big Ulrich vector bundles on surfaces and threefolds}.
Internat. J. Math. \textbf{32} (2021), no.~14, Paper No.~2150111, 18 pp.

\bibitem[LMS1]{lms1} A.~F.~Lopez, R.~Mu\~{n}oz, J.~C.~Sierra. 
\textit{On the classification of non-big Ulrich vector bundles on fourfolds}.
Ann. Sc. Norm. Super. Pisa Cl. Sci. (5) \textbf{26} (2025), no.~2, 707-755.

\bibitem[LMS2]{lms2} A.~F.~Lopez, R. Mu\~{n}oz, J.C. Sierra.
\textit{Non-big Ulrich bundles: the classification on quadrics and the case of small numerical dimension}.
Manuscripta Math. \textbf{174} (2024), 517-533.

\bibitem[LP]{lp} A.~Lanteri, M.~Palleschi.
\textit{Projective manifolds containing many rational curves}. 
Indiana Univ. Math. J. \textbf{36} (1987), no.~4, 857-865. 

\bibitem[LR]{lr} A.~F.~Lopez, D.~Raychaudhury.
\textit{Ulrich subvarieties and the non-existence of low rank Ulrich bundles on complete intersections}.
Trans. Amer. Math. Soc. \textbf{379} (2026), no.\ 8, 5809-5846.

\bibitem[LRT]{lrt} A.~F.~Lopez, D.~Raychaudhury, Y. Takahashi.
\textit{A linear lower bound on the Ulrich complexity of hypersurfaces}.
Preprint 2026, arXiv: 2607.03944.

\bibitem[LS]{ls} A.~F.~Lopez, J.~C.~Sierra. 
\textit{A geometrical view of Ulrich vector bundles}.
Int. Math. Res. Not. IMRN(2023), no.~11, 9754-9776.

\bibitem[MOS]{mos} R.~Mu\~{n}oz, G.~Occhetta, L.~E.~Sol\'a Conde.
\textit{On rank 2 vector bundles on Fano manifolds}. 
Kyoto J. Math. \textbf{54} (2014), no. 1, 167-197.

\bibitem[PZ]{pz} Y.~Prokhorov, M.~Zaidenberg. 
\textit{Examples of cylindrical Fano fourfolds}.
Eur. J. Math. \textbf{2} (2016), no.~1, 262-282.

\bibitem[R]{ru} F.~Russo.
\textit{On the geometry of some special projective varieties}. 
Lecture Notes of the Unione Matematica Italiana, \textbf{18}. Springer, Cham; Unione Matematica Italiana, Bologna, 2016.

\bibitem[ST]{st} 
\textit{The Stacks project}. 
\url{https://stacks.math.columbia.edu}

\bibitem[W]{wa} K.~Watanabe.
\textit{$\P^1$-bundles admitting another smooth morphism of relative dimension one}. 
J. Algebra \textbf{414} (2014), 105-119.
 
\end{thebibliography}
\end{document}